\documentclass[11pt]{article} 
\usepackage{latexsym} 
\newcommand{\eproof}{\mbox{\ }\hfill $\Box$ \par \vskip 10pt}

\newtheorem{Theorem}{Theorem}[section] 
\newtheorem{lemma}[Theorem]{Lemma} 
\newtheorem{prop}[Theorem]{Proposition}

\begin{document}

\title{Dispersive estimates of solutions to the
wave equation with a potential in dimensions $n\ge 4$}

\author{{\sc Georgi Vodev}}

\date{} 
 
\maketitle 

\begin{abstract}
\noindent
We prove dispersive estimates 
for solutions to the wave  
equation with a real-valued potential
$V\in L^\infty({\bf R}^n)$, $n\ge 4$, satisfying 
$V(x)=O(\langle x\rangle^{-(n+1)/2-\epsilon})$, $\epsilon>0$. 
\end{abstract}

\setcounter{section}{0}
\section{Introduction and statement of results}

Let $V\in L^\infty({\bf R}^n)$, $n\ge 4$,
 be a real-valued function satisfying
$$\left|V(x)\right|\le C\langle x\rangle^{-\delta},
\quad \forall x\in {\bf R}^n,\eqno{(1.1)}$$
with constants $C>0$ and $\delta>(n+1)/2$, 
where $\langle x\rangle=(1+|x|^2)^{1/2}$. 
Denote by $G_0$ and $G$ the self-adjoint realizations of the
operators $-\Delta$ and $-\Delta+V(x)$ on $L^2({\bf R}^n)$.
It is well known that the absolutely continuous spectrums of the
operators $G_0$ and $G$ coincide with the interval $[0,+\infty)$.
Moreover, by Kato's theorem the operator $G$ has no strictly 
positive eigenvalues. This implies that  
$G$ has no strictly positive resonances neither (e.g. see 
\cite{kn:V}). 
Throughout this paper, given $1\le p\le +\infty$, $L^p$ will 
denote the space $L^p({\bf R}^n)$. Also, 
given an $a>0$ denote by $\chi_a\in C^\infty({\bf R})$ 
 a real-valued function supported in the interval
$[a,+\infty)$, $\chi_a=1$ on $[2a,+\infty)$. 
Our main result is the following

\begin{Theorem}
Assume (1.1) fulfilled. Then,  
for every $a>0$, $2\le p<\frac{2(n-1)}{n-3}$, 
there exists a constant $C>0$ so that the following estimate holds
$$\left\|e^{it\sqrt{G}}(\sqrt{G})^{-\alpha(n+1)/2}
\chi_a(\sqrt{G})\right\|_{L^{p'}\to L^p}
\le C|t|^{-\alpha(n-1)/2},\quad \forall t\neq 0,\eqno{(1.2)}$$
where $1/p+1/p'=1$, $\alpha=1-2/p$. Moreover, 
for every $a>0$, $2\le p<+\infty$, we have
$$\left\|e^{it\sqrt{G}}(\sqrt{G})^{-\alpha(n-1)}
\chi_a(\sqrt{G})\right\|_{L^{p'}\to L^p}
\le C|t|^{-\alpha(n-1)/2},\quad \forall t\neq 0,\eqno{(1.3)}$$
where $1/p+1/p'=1$, $\alpha=1-2/p$, and, $\forall 0<\epsilon\ll 1$, 
$$\left\|e^{it\sqrt{G}}(\sqrt{G})^{-\alpha(n+1)/2}
\chi_a(\sqrt{G})\langle x\rangle^{-\alpha(n/2+\epsilon)}
\right\|_{L^2\to L^p}\le C_\epsilon|t|^{-\alpha(n-1)/2},
\quad \forall t\neq 0.\eqno{(1.4)}$$
\end{Theorem}

{\bf Remark 1.} It is easy to see from the proof that the 
estimates (1.3) and (1.4) hold true for $p=+\infty$, $p'=1$, with 
$(\sqrt{G})^{-\alpha}$ replaced by $(\sqrt{G})^{-1-\epsilon}$, 
$\forall\epsilon>0$. 

{\bf Remark 2.} It is also clear from the proof that the
following intermediate estimate between (1.2) and (1.3) holds:
$$\left\|e^{it\sqrt{G}}(\sqrt{G})^{-\alpha((n+1)/2+q)}
\chi_a(\sqrt{G})\right\|_{L^{p'}\to L^p}
\le C|t|^{-\alpha(n-1)/2},\quad \forall t\neq 0,\eqno{(1.5)}$$
for every $0<q<(n-3)/2$ and for $2\le p<\frac{2(n-1-2q)}{n-3-2q}$.

{\bf Remark 3.} The desired result would be to prove (1.2) for
$n\ge 4$ and for all $2\le p<+\infty$, but we doubt that such a
statement can be proved without imposing extra assumptions on the potential.
Since in (1.3) we have, roughly speaking, a loss of $(n-3)/2$ derivatives, 
it is natural to expect that it suffices to control the behaviour of
$(n-3)/2$ (radial) derivatives of $V$ in order that (1.2) holds true 
for all $2\le p<+\infty$. In the case of the Schr\"odinger group Goldberg
and Visan \cite{kn:GoV} have recently showed that there exist compactly
supported potentials $V\in C^k({\rm R}^n)$, $\forall k<(n-3)/2$,
for which the optimal (without loss of derivatives) $L^1\to L^\infty$
dispersive estimate fails to hold, and this is a typical high frequency
phenomenon. Similar thing should occur in the case of the wave group, too. 

{\bf Remark 4.} Note that given a smooth, bounded function 
$f$ supported in the interval $(0,+\infty)$,
the operator-valued function $f(\sqrt{G})$ is well defined even if
the operator $G$ is not non-negative. In particular, the operators 
above are well defined. 

{\bf Remark 5.} The operator $G$ may have in general a finite number
of eigenvalues $\lambda_j\le 0$ as well as a resonance at zero.
We eliminate the influence of these points to the decay properties
of the wave group by cutting off with the operator
$\chi_a(\sqrt{G})$. It is natural to expect that if the zero is
neither an eigenvalue nor a resonance of $G$, the statements
of Theorem 1.1 hold true with $\chi_a$ replaced by the
characteristic function, $\chi$, of the interval $[0,+\infty)$
(the absolutely continuous spectrum of $G$) and 
$(\sqrt{G})^{-(n-1)}$ replaced by $(\sqrt{G})^{-(n+1)/2}\left\langle G
\right\rangle^{-(n-3)/4}$. To prove this, it suffices to show
that the estimate (1.2) holds in this case for all
$2\le p<+\infty$ with $\chi_a$ replaced
by $\chi(1-\chi_a)$ with some $a>0$ small enough. The proof of 
such an estimate, however, requires different techniques than those
developed in the present paper. 

It is well known that (1.2) for all $2\le p<+\infty$ holds true for 
the free operator $G_0$ with $\chi_a\equiv 1$ in all dimensions. 
Recently, (1.2) for all $2\le p<+\infty$ has been proved in \cite{kn:CCV}
when $n=2$ and $n=3$ 
by a different method (using some properties of the resolvent of 
$G_0$ which are no longer valid when $n\ge 4$). 
For $n=3$, an analogue of (1.2) 
(for $2\le p\le 4$) is proved by Georgiev and Visciglia \cite{kn:GV}
for non-negative potentials satisfying (1.1) as well as an extra
regularity assumption. Beals and Strauss \cite{kn:BS} proved
 an analogue of (1.2) in all dimensions $n\ge 3$ 
(for $2\le p\le \frac{2(n+1)}{n-1}$) for
a class of non-negative potentials decaying much faster at infinity 
than those we consider in the present paper, while Beals \cite{kn:B} 
proved the estimate (1.2) with a loss of $\varepsilon$-derivatives  
for $n\ge 3$ and all $2\le p\le +\infty$ for potentials belonging to the
Schwartz class ${\cal S}({\bf R}^n)$. Recently, in
\cite{kn:CuV} an analogue of (1.2) with $n=3$ has been proved for a
class of potentials satisfying (1.1) with $4/3<\delta\le 2$, but 
with a weaker decay as $|t|\to +\infty$. Note also the work \cite{kn:CV},
where a better time decay than that in (1.2) has been obtained on
weighted $L^p$ spaces for potentials satisfying (1.1) with $n=3$
and $\delta>2$. In the case $n=3$ D'ancona and Pierfelice
\cite{kn:DP} obtained dispersive estimates for real-valued potentials, 
$V$, belonging to a quite large subset of the Kato class (with a 
small Kato norm of the negative part of $V$). 
Assuming additionally that $G$ has no real
resonances, they proved the estimate
$$\left\|\frac{\sin\left(t\sqrt{G}\right)}{\sqrt{G}}\right
\|_{\dot B^1_{1,1}\to L^\infty}
\le C|t|^{-1},\quad \forall t\neq 0,\eqno{(1.6)}$$
where $\dot B^1_{1,1}$ is the homogeneous Besov space. Note that a similar
estimate to (1.6) ( with $\dot B^1_{1,1}$ replaced by $W^{1,1}$)
has been recently proved by Krieger and Schlag (see Section 5 of
\cite{kn:KS}) still in dimension $n=3$ for potentials satisfying
(1.1) with $\delta>3$.

Our method is quite different from those used in the papers mentioned
above. It consists of reducing the estimates (1.2)-(1.4) to
{\it semi-classical} ones (see Theorem 4.1) valid for all
$2\le p\le +\infty$, which in turn makes easier applying interpolation
arguments. On the other hand, the estimates of Theorem 4.1 are
reduced to uniform $L^2\to L^2$ estimates for the operator
$e^{it\sqrt{G}}\varphi(h\sqrt{G})$, $0<h\le 1$, $\varphi\in
C_0^\infty((0,+\infty))$. These latter estimates are proved 
in Section 3, using the properties of both the free and the perturbed 
resolvents on weighted $L^2$ spaces. In Section 2 we prove some properties 
of the free wave group $e^{it\sqrt{G_0}}$ (see Proposition 2.1),
which play an essential role in our proof of the above estimates.
Note that the method presented here works also for $n=2$ and $n=3$, but
in this case a simpler proof (of (1.2) for all 
$2\le p<+\infty$) is given in \cite{kn:CCV}. That is why in the present
paper we will treat the case of $n\ge 4$ only.

Similar ideas have been used in \cite{kn:V2} to prove 
$L^1\to L^\infty$ dispersive estimates for the Schr\"odinger group
$e^{itG}$ in the case $n=3$ for potentials satisfying (1.1) with
$\delta>5/2$. The method of \cite{kn:V2} actually works in all 
dimensions $n\ge 3$ for potentials satisfying (1.1) with
$\delta>(n+2)/2$ (see \cite{kn:V3})
 and gives $L^1\to L^\infty$ dispersive estimates for 
$e^{itG}$ with the right time decay for $|t|\gg 1$, but with a loss
of $(n-3)/2$ derivatives similarly to the dispersive estimate
(1.3) for the wave group.

\section{Preliminary estimates}

The following properties of the free wave group 
will play a key role in the proof of our dispersive estimates.

\begin{prop}
Let $\varphi\in C_0^\infty((0,+\infty))$. For every 
$0\le s\le (n-1)/2$, $0<\epsilon\ll 1$, $0<h\le 1$ , we have 
$$\left\|\langle x\rangle^{-s}
e^{it\sqrt{G_0}}\varphi(h\sqrt{G_0})\langle x\rangle^{-s}
\right\|_{L^2\to L^2}\le C\langle t\rangle^{-s},\quad 
\forall t,\eqno{(2.1)}$$
$$\left\|e^{it\sqrt{G_0}}\varphi(h\sqrt{G_0})
\right\|_{L^1\to L^\infty}\le
 Ch^{-(n+1)/2}|t|^{-(n-1)/2},\quad \forall t\neq 0,\eqno{(2.2)}$$
$$\left\|e^{it\sqrt{G_0}}\varphi(h\sqrt{G_0})
\langle x\rangle^{-1/2-s-\epsilon}\right\|_{L^2\to L^\infty}\le
 Ch^{-(n+1)/2}|t|^{-s},\quad \forall t\neq 0,\eqno{(2.3)}$$
$$\int_{-\infty}^\infty|t|^{2s}\left\|\langle x\rangle^{-1/2-s-\epsilon}
e^{it\sqrt{G_0}}\varphi(h\sqrt{G_0})f\right\|^2_{L^2}dt
\le Ch^{-n}\|f\|^2_{L^1},\quad \forall f\in L^1,\eqno{(2.4)}$$ 
with a constant $C>0$ independent of $t$, $h$ and $f$.
\end{prop}

{\it Proof.} We will prove these estimates for all $n\ge 2$. 
We are going to take advantage of the formula
$$e^{it\sqrt{G_0}}\varphi(h\sqrt{G_0})=(\pi i)^{-1}\int_0^\infty
 e^{it\lambda}\varphi(h\lambda)\left(R_0^+(\lambda)-R_0^-(\lambda)\right)
\lambda d\lambda,\eqno{(2.5)}$$
where $R_0^\pm(\lambda)$ are the outgoing and incoming free resolvents
with kernels given in terms of the Hankel functions by
$$[R^\pm_0(\lambda)](x,y)=\pm\frac{i}{4}\left(\frac{\lambda}
{2\pi |x-y|}\right)^\nu H_\nu^\pm(\lambda|x-y|),$$
where $\nu=(n-2)/2$. Hence, the kernel of the operator
$e^{it\sqrt{G_0}}\varphi(h\sqrt{G_0})$ is of the form
$K_h(|x-y|,t)$, where
$$K_h(\sigma,t)=\frac{\sigma^{-2\nu}}{(2\pi)^{\nu+1}}
\int_0^\infty
 e^{it\lambda}\varphi(h\lambda){\cal J}_\nu(\sigma\lambda)
\lambda d\lambda,\eqno{(2.6)}$$
where ${\cal J}_\nu(z)=z^\nu J_\nu(z)$, $J_\nu(z)=(H^+_\nu(z)+
H^-_\nu(z))/2$ is the Bessel function of order $\nu$. 
We will first show that the above estimates follow from the following

\begin{lemma} For every $0\le s\le (n-1)/2$, $\sigma>0$, $t\neq 0$,
$0<h\le 1$, we have
$$\left|K_h(\sigma,t)\right|\le C|t|^{-s}h^{-(n+1)/2}\sigma^{-(n-1)/2+s},
\eqno{(2.7)}$$
 $$\int_{-\infty}^\infty|t|^{2s}\left|K_h(\sigma,t)\right|^2dt
\le Ch^{-n}\sigma^{-(n-1)+2s},\eqno{(2.8)}$$
 $$\int_0^{|t|/2}\sigma^{n-1}\left|K_h(\sigma,t)\right|d\sigma
\le C|t|^{-(n-1)/2},\quad|t|\ge 1,\eqno{(2.9)}$$
with a constant $C>0$ independent of $\sigma$, $t$ and $h$.
\end{lemma}

Clearly, (2.1) is trivial for $|t|\le 1$, so we may suppose that
$|t|\ge 1$. Given a set ${\cal M}\subset{\bf R}^n$ denote by
$\eta({\cal M})$ the characteristic function of ${\cal M}$. We have
$$\left\|\langle x\rangle^{-s}
e^{it\sqrt{G_0}}\varphi(h\sqrt{G_0})\langle x\rangle^{-s}
\right\|_{L^2\to L^2}$$ $$\le \left\|\eta(|x|\le |t|/4)
e^{it\sqrt{G_0}}\varphi(h\sqrt{G_0})\eta(|x|\le |t|/4)
\right\|_{L^2\to L^2}+C\langle t\rangle^{-s}.\eqno{(2.10)}$$
In view of Schur's lemma the norm in the RHS of (2.10) is upper bounded
by
$$\sup_{|x|\le|t|/4}\int_{|y|\le|t|/4}|K_h(|x-y|,t)|dy\le
\int_{|\xi|\le|t|/2}|K_h(|\xi|,t)|d\xi$$ $$\le
 C\int_0^{|t|/2}\sigma^{n-1}\left|K_h(\sigma,t)\right|d\sigma
\le C|t|^{-(n-1)/2},\eqno{(2.11)}$$
where we have used (2.9). Now (2.1) follows from (2.10) and (2.11).
Clearly, (2.2) follows from (2.7) with $s=(n-1)/2$. 
To prove (2.3) observe that in view of (2.7) we have
$$\left\|e^{it\sqrt{G_0}}\varphi(h\sqrt{G_0})
\langle x\rangle^{-1/2-s-\epsilon}\right\|^2_{L^2\to L^\infty}\le
\sup_{y\in{\bf R}^n}\int_{{\bf R}^n}|K_h(|x-y|,t)|^2
\langle x\rangle^{-1-2s-2\epsilon}dx$$
 $$\le Ch^{-n-1}|t|^{-2s}\sup_{y\in{\bf R}^n}\int_{{\bf R}^n}
|x-y|^{-(n-1)+2s}\langle x\rangle^{-1-2s-2\epsilon}dx\le 
 Ch^{-n-1}|t|^{-2s},\eqno{(2.12)}$$
which is the desired result. To prove (2.4) observe that for any 
$f\in L^1$, we have
$$\left|\langle x\rangle^{-1/2-s-\epsilon}e^{it\sqrt{G_0}}
\varphi(h\sqrt{G_0})f\right|\le \int_{{\bf R}^n}
\langle x\rangle^{-1/2-s-\epsilon}|K_h(|x-y|,t)||f(y)|dy$$
 $$\le \left(\int_{{\bf R}^n}
\langle x\rangle^{-1-2s-2\epsilon}|K_h(|x-y|,t)|^2|f(y)|
dy\right)^{1/2}\left(\int_{{\bf R}^n}|f(y)|dy\right)^{1/2}.$$
Hence,
$$\left\|\langle x\rangle^{-1/2-s-\epsilon}
e^{it\sqrt{G_0}}\varphi(h\sqrt{G_0})f\right\|^2_{L^2}\le
\|f\|_{L^1}\int_{{\bf R}^n}\int_{{\bf R}^n}
\langle x\rangle^{-1-2s-2\epsilon}|K_h(|x-y|,t)|^2|f(y)|dx dy.$$
Thus we obtain
$$\int_{-\infty}^\infty|t|^{2s}\left\|\langle x\rangle^{-1/2-s-\epsilon}
e^{it\sqrt{G_0}}\varphi(h\sqrt{G_0})f\right\|^2_{L^2}dt$$ $$\le 
\|f\|_{L^1}\int_{{\bf R}^n}\int_{{\bf R}^n}
\langle x\rangle^{-1-2s-2\epsilon}\int_{-\infty}^\infty|t|^{2s}
|K_h(|x-y|,t)|^2dt|f(y)|dx dy$$ 
 $$\le\|f\|^2_{L^1}\sup_{y\in{\bf R}^n}\int_{{\bf R}^n}
\langle x\rangle^{-1-2s-2\epsilon}\int_{-\infty}^\infty|t|^{2s}
|K_h(|x-y|,t)|^2dtdx$$ 
 $$\le Ch^{-n}\|f\|^2_{L^1}\sup_{y\in{\bf R}^n}\int_{{\bf R}^n}
\langle x\rangle^{-1-2s-2\epsilon}|x-y|^{-(n-1)+2s}dx\le 
 Ch^{-n}\|f\|^2_{L^1},\eqno{(2.13)}$$
where we have used (2.8), which is the desired result.\\

{\it Proof of Lemma 2.2.} It is well known that the function 
${\cal J}_\nu(z)$ satisfies the bounds
$$\left|\partial_z^k{\cal J}_\nu(z)\right|\le C_kz^{(n-3)/2},\quad
 z\ge 1,\eqno{(2.14)}$$
for every integer $k\ge 0$, while near $z=0$ the function
${\cal J}_\nu(z)$ is equal to $z^{2\nu}$ times an even analytic function.
Therefore, when $n\ge 3$, we have
$$\left|\partial_z^k{\cal J}_\nu(z)\right|\le Cz^{n-2-k},\quad
 0<z\le 1,\eqno{(2.15)}$$
for all integers $0\le k\le n-2$. For $n=2$, we have
$$\left|\partial_z^k{\cal J}_0(z)\right|\le Cz^{k},\quad
 0<z\le 1,\,k=0,1.\eqno{(2.16)}$$
Thus, when $n\ge 3$, we obtain
$$\left|\partial_z^k{\cal J}_\nu(z)\right|\le Cz^{n-2-k}\langle z
\rangle^{-(n-1)/2+k},\quad\forall z>0,\eqno{(2.17)}$$
for all integers $0\le k\le n-2$. For $n=2$, we obtain
$$\left|\partial_z^k{\cal J}_0(z)\right|\le Cz^{k}
\langle z\rangle^{-1/2-k},\quad\forall z>0,\,k=0,1.\eqno{(2.18)}$$
Let $n\ge 3$ and let $m$ be any integer such that $0\le m\le n-2$.
Integrating by parts $m$ times the integral in (2.6) and using (2.17),
we obtain
$$\left|t^mK_h(\sigma,t)\right|\le h^{-1}\sigma^{-2\nu}\left|\int_0^\infty
 e^{it\lambda}\partial^m_\lambda\left(\widetilde\varphi(h\lambda){\cal J}_\nu
(\sigma\lambda)\right)d\lambda\right|$$
 $$\le Ch^{-1}\sigma^{-2\nu}\sum_{k=0}^mh^{m-k}\sigma^k\int_0^\infty
\left|\frac{\partial^{m-k}\widetilde\varphi}{\partial\lambda^{m-k}}(h\lambda)
\right| \left|\frac{\partial^{k}{\cal J}_\nu}{\partial\lambda^{k}}
(\sigma\lambda)\right|d\lambda$$
  $$\le Ch^{-2}\sigma^{-2\nu}\sum_{k=0}^mh^{m-k}\sigma^k\int_{{\rm supp}\,
\varphi}\left|\frac{\partial^{k}{\cal J}_\nu}{\partial\lambda^{k}}
(\sigma\lambda/h)\right|d\lambda$$
 $$\le Ch^{-2}\sigma^{-2\nu}\sum_{k=0}^mh^{m-k}\sigma^k(\sigma/h)^{n-2-k}
\langle \sigma/h\rangle^{-(n-1)/2+k}$$ $$\le Ch^{m-n}
\langle \sigma/h\rangle^{-(n-1)/2+m},\eqno{(2.19)}$$
where $\widetilde\varphi(\lambda)=\lambda\varphi(\lambda)$. Clearly, 
(2.19) holds for all real $0\le m\le n-2$, and in particular
for $m=s$, $0\le s\le (n-1)/2$. When $n=2$ one can see in the same way,
using (2.18) instead of (2.17), that (2.19) still holds for $m=0$ and
$m=1$, and hence for all real $0\le m\le 1$. In particular, (2.19)
holds in this case with $m=s$, $0\le s\le 1/2$. Now (2.7) follows from
(2.19) (with $m=s$) and the inequality
$$\langle \sigma/h\rangle^{-(n-1)/2+s}\le h^{(n-1)/2-s}
\sigma^{-(n-1)/2+s},\quad 0\le s\le (n-1)/2.$$
Let $n\ge 3$ and let $m$ be any integer such that $0\le m\le n-2$.
By Plancherel's identity and (2.6), (2.17), we obtain
$$\int_{-\infty}^\infty|t|^{2m}\left|K_h(\sigma,t)\right|^2dt=
Ch^{-2}\sigma^{-4\nu}\int_0^\infty\left|
\partial^m_\lambda\left(\widetilde\varphi(h\lambda){\cal J}_\nu
(\sigma\lambda)\right)\right|^2d\lambda$$
 $$\le Ch^{-2}\sigma^{-4\nu}\sum_{k=0}^mh^{2(m-k)}\sigma^{2k}\int_0^\infty
\left|\frac{\partial^{m-k}\widetilde\varphi}{\partial\lambda^{m-k}}(h\lambda)
\right|^2 \left|\frac{\partial^{k}{\cal J}_\nu}{\partial\lambda^{k}}
(\sigma\lambda)\right|^2d\lambda$$
  $$\le Ch^{-3}\sigma^{-4\nu}\sum_{k=0}^mh^{2(m-k)}\sigma^{2k}
\int_{{\rm supp}\,
\varphi}\left|\frac{\partial^{k}{\cal J}_\nu}{\partial\lambda^{k}}
(\sigma\lambda/h)\right|^2d\lambda$$
 $$\le Ch^{-3}\sigma^{-4\nu}\sum_{k=0}^mh^{2(m-k)}\sigma^{2k}
(\sigma/h)^{2(n-2-k)}
\langle \sigma/h\rangle^{-(n-1)+2k}$$ $$\le Ch^{2m-2n+1}
\langle \sigma/h\rangle^{-(n-1)+2m}.\eqno{(2.20)}$$
Clearly, (2.20) holds for all real $0\le m\le n-2$, and in particular
for $m=s$, $0\le s\le (n-1)/2$. When $n=2$ it is easy to see, using (2.18)
instead of (2.17), that (2.20) still holds for $m=0$ and $m=1$,
and hence for all real $0\le m\le 1$. Now (2.8) follows from (2.20)
in the same way as above. 
To prove (2.9) we will use that ${\cal J}_\nu(z)=e^{iz}b_\nu^++
e^{-iz}b_\nu^-$, where $b_\nu^\pm(z)$ are symbols of order $(n-3)/2$
for $z\gg 1$, and satisfy the bounds
$$\left|\partial_z^kb_\nu^\pm(z)\right|\le C_kz^{(n-3)/2-k},\quad
 z\ge z_0,\eqno{(2.21)}$$
for every integer $k\ge 0$ and every $z_0>0$ with a constant $C_k>0$
depending on $z_0$. Write the function $K_h$ in the form
$K_h^++K_h^-$, where
$$K_h^\pm(\sigma,t)=\frac{h^{-1}\sigma^{-2\nu}}{(2\pi)^{\nu+1}}
\int_0^\infty
 e^{i(t\pm\sigma)\lambda}\widetilde\varphi(h\lambda)b_\nu^\pm(\sigma\lambda)
 d\lambda,\eqno{(2.22)}$$
$\widetilde\varphi$ being as above. By (2.7) we have
$$\int_0^h\sigma^{n-1}|K_h(\sigma;t)|d\sigma
\le Ch^{(n-1)/2}|t|^{-(n-1)/2},\eqno{(2.23)}$$
so it sufices to show that
$$\int_h^{|t|/2}\sigma^{n-1}|K_h^\pm(\sigma;t)|d\sigma
\le C|t|^{-(n-1)/2}.\eqno{(2.24)}$$
Integrating  by parts $n$ times the integral in (2.22) and using (2.21),
we obtain, for $h\le\sigma\le |t|/2$,
$$\left|t^nK_h^\pm(\sigma,t)\right|\le Ch^{-1}\sigma^{-2\nu}
\int_0^\infty\left|\partial^n_\lambda\left(\widetilde
\varphi(h\lambda)b^\pm_\nu
(\sigma\lambda)\right)\right|d\lambda$$
 $$\le Ch^{-1}\sigma^{-2\nu}\sum_{k=0}^nh^{n-k}\sigma^k\int_0^\infty
\left|\frac{\partial^{n-k}\widetilde\varphi}{\partial\lambda^{n-k}}(h\lambda)
\right| \left|\frac{\partial^{k}b^\pm_\nu}{\partial\lambda^{k}}
(\sigma\lambda)\right|d\lambda$$
 $$\le Ch^{-2}\sigma^{-2\nu}\sum_{k=0}^nh^{n-k}\sigma^k\int_{{\rm supp}\,
\varphi}\left|\frac{\partial^{k}b^\pm_\nu}{\partial\lambda^{k}}
(\sigma\lambda/h)\right|d\lambda$$
 $$\le Ch^{-2}\sigma^{-2\nu}\sum_{k=0}^nh^{n-k}
\sigma^k(\sigma/h)^{(n-3)/2-k}\le Ch^{(n-1)/2}
\sigma^{-(n-1)/2}.\eqno{(2.25)}$$
Hence
$$\int_h^{|t|/2}\sigma^{n-1}|K_h^\pm(\sigma;t)|d\sigma
\le Ch^{(n-1)/2}|t|^{-n}\int_h^{|t|/2}\sigma^{(n-1)/2}d\sigma
\le Ch^{(n-1)/2}|t|^{-(n-1)/2},$$
which implies (2.24).
\eproof

We will also need the following

\begin{lemma} Assume (1.1) fulfilled. Then, for every 
$\varphi\in C_0^\infty((0,+\infty))$, $0\le s\le \delta$, 
$1\le p\le \infty$, $0<h\le 1$, we have
$$\left\|\langle x\rangle^{-s}\varphi(h\sqrt{G_0})
\langle x\rangle^{s}\right\|_{L^2\to L^2}\le C,\eqno{(2.26)}$$
 $$\left\|\langle x\rangle^{-s}\varphi(h\sqrt{G})
\langle x\rangle^{s}\right\|_{L^2\to L^2}\le C,\eqno{(2.27)}$$
 $$\left\|\left(\varphi(h\sqrt{G_0})
-\varphi(h\sqrt{G})\right)
\langle x\rangle^{s}\right\|_{L^2\to L^2}\le Ch^2,\eqno{(2.28)}$$
 $$\left\|\varphi(h\sqrt{G_0})\right\|_{L^p\to L^p}\le C,\eqno{(2.29)}$$
 $$\left\|\varphi(h\sqrt{G})\right\|_{L^p\to L^p}\le C,\eqno{(2.30)}$$
 $$\left\|\varphi(h\sqrt{G})-\varphi(h\sqrt{G_0})
\right\|_{L^p\to L^p}\le Ch^2,\eqno{(2.31)}$$
 $$\left\|\varphi(h\sqrt{G_0})\right\|_{L^2\to L^p}\le Ch^{-n|\frac{1}{2}
-\frac{1}{p}|},\eqno{(2.32)}$$
 $$\left\|\varphi(h\sqrt{G})\right\|_{L^2\to L^p}\le C
h^{-n|\frac{1}{2}-\frac{1}{p}|},\eqno{(2.33)}$$
 $$\left\|\left(\varphi(h\sqrt{G})-\varphi(h\sqrt{G_0})\right)
\langle x\rangle^{s}\right\|_{L^2\to L^p}\le C
h^{2-n|\frac{1}{2}-\frac{1}{p}|},\eqno{(2.34)}$$
with a constant $C>0$ independent of $h$.
\end{lemma}

{\it Proof.} Define the function $\psi\in C_0^\infty
((0,+\infty))$ by $\psi(\sigma^2)=\varphi(\sigma)$. We
will make use of the Helffer-Sj\"ostrand formula
$$\psi(h^2G)=\frac{1}{\pi}\int_{\bf C}\frac{\partial\widetilde\psi}
{\partial\overline z}(z)(h^2G-z)^{-1}L(dz),\eqno{(2.35)}$$
where $L(dz)$ denotes the Lebesgue measure on ${\bf C}$, and
$\widetilde\psi\in C_0^\infty({\bf C})$ 
is an almost analytic continuation of
$\psi$ supported in a small complex neighbourhood of supp$\,\psi$
and satisfying
$$\left|\frac{\partial\widetilde\psi}
{\partial\overline z}(z)\right|
\le C_N|{\rm Im}\,z|^N,\quad\forall N\ge 1.$$
Therefore, (2.26) would follow from the bound
$$\left\|\langle x\rangle^{-s}(h^2G_0-z)^{-1}
\langle x\rangle^s\right\|_{L^2\to L^2}
\le C|{\rm Im}\,z|^{-q_1},\eqno{(2.36)}$$
for $z\in{\bf C}_\psi:={\rm supp}\,\widetilde\psi$, 
${\rm Im}\,z\neq 0$, with some constants 
$C,q_1>0$ independent of $h$ and $z$.
To prove (2.36) we will use the identity
$$\langle x\rangle^{-s}(h^2G_0-z)^{-1}\langle x\rangle^s=
(h^2G_0-z)^{-1}+h^2\langle x\rangle^{-s}(h^2G_0-z)^{-1}
[\Delta,\langle x\rangle^s](h^2G_0-z)^{-1},$$
together with the bound
$$\left\|(h\nabla_x)^j(h^2G_0-z)^{-1}\right\|_{L^2\to L^2}
\le C|{\rm Im}\,z|^{-1},\quad j=0,1.\eqno{(2.37)}$$
Thus we obtain
$$\left\|\langle x\rangle^{-s}
(h^2G_0-z)^{-1}\langle x\rangle^s\right\|_{L^2\to L^2}
\le C\left\|(h^2G_0-z)^{-1}\right\|_{L^2\to L^2}$$ 
$$ +Ch^2\left\|\langle x\rangle^{-s}
(h^2G_0-z)^{-1}\langle x\rangle^{s-1}\right\|_{L^2\to L^2}
\left(\left\|\nabla_x(h^2G_0-z)^{-1}
\right\|_{L^2\to L^2}+\left\|(h^2G_0-z)^{-1}\right\|_{L^2\to L^2}\right)$$
 $$\le C|{\rm Im}\,z|^{-1}\left(1+h\left\|\langle x\rangle^{-s}
(h^2G_0-z)^{-1}\langle x\rangle^{s-1}\right\|_{L^2\to L^2}\right).$$
Repeating this a finite number of times leads to (2.36).

The estimate (2.27) follows from (2.26) and (2.28). To prove (2.28)
we use (2.35) to obtain
$$\left\|\left(\psi(h^2G)-
\psi(h^2G_0)\right)\langle x\rangle^s\right\|_{L^2\to L^2}$$ 
  $$\le O(h^2)\int_{\bf C}\left|\frac{\partial\widetilde\psi}
{\partial\overline z}(z)\right|\left\|
(h^2G-z)^{-1}V(h^2G_0-z)^{-1}\langle x\rangle^s
\right\|_{L^2\to L^2}L(dz)$$
 $$\le O_N(h^2)\int_{{\bf C}_\psi}\left|{\rm Im}\,z
\right|^N\left\|
(h^2G-z)^{-1}\right\|_{L^2\to L^2}\left\|\langle x\rangle^{-s}
(h^2G-z)^{-1}\langle x\rangle^s\right\|_{L^2\to L^2}L(dz),$$
which together with
(2.36) and (2.43) below imply the desired result.

The estimate (2.29) follows from (2.35) and the bound
$$\left\|(h^2G_0-z)^{-1}\right\|_{L^p\to L^p}
\le C|{\rm Im}\,z|^{-q_2},\eqno{(2.38)}$$
for $z\in{\bf C}_\psi$, ${\rm Im}\,z\neq 0$, with some constants 
$C,q_2>0$ independent of $h$ and $z$. 
The estimate (2.30) follows from (2.29) and (2.31). 
Using (2.35) as above, we obtain
$$\left\|\psi(h^2G)-
\psi(h^2G_0)\right\|_{L^p\to L^p}$$ $$
  \le O_N(h^2)\int_{{\bf C}_\psi}\left|{\rm Im}\,z
\right|^N\left\|
(h^2G-z)^{-1}\right\|_{L^p\to L^p}\left\|
(h^2G_0-z)^{-1}\right\|_{L^p\to L^p}L(dz).\eqno{(2.39)}$$
Thus, (2.31) follows from (2.38), (2.39) and the bound
$$\left\|(h^2G-z)^{-1}\right\|_{L^p\to L^p}
\le C|{\rm Im}\,z|^{-q_3},\eqno{(2.40)}$$
for $z\in{\bf C}_\psi$, ${\rm Im}\,z\neq 0$, with some constants 
$C,q_3>0$ independent of $h$ and $z$. 
Observe next that by the resolvent identity we can write
$$(h^2G-z)^{-1}=\sum_{j=0}^M h^{2j}\left((h^2G_0-z)^{-1}V\right)^j
(h^2G_0-z)^{-1}$$ $$+h^{2M+2}\left((h^2G_0-z)^{-1}V\right)^M(h^2G-z)^{-1}V
(h^2G_0-z)^{-1}\eqno{(2.41)}$$
for every integer $M\ge 1$. Taking $M$ big enough, it is easy to
see that (2.40) follows from (2.41) together with (2.38) and the
following well known bounds
$$\left\|(h^2G_0-z)^{-1}\right\|_{L^p\to L^2}+
\left\|(h^2G_0-z)^{-1}\right\|_{L^2\to L^p}
\le Ch^{-q}|{\rm Im}\,z|^{-q-1},\eqno{(2.42)}$$
for $z\in{\bf C}_\psi$, ${\rm Im}\,z\neq 0$, with a constant 
$C>0$ independent of $h$ and $z$, where $q=n\left|\frac{1}{2}-
\frac{1}{p}\right|$, and 
$$\left\|(h^2G-z)^{-1}\right\|_{L^2\to L^2}
\le |{\rm Im}\,z|^{-1}.\eqno{(2.43)}$$
The estimate (2.32) follows from (2.35) and (2.42), while 
(2.33) follows from (2.32) and (2.34). 
To prove (2.34), we use (2.35) to obtain
$$\left\|\left(\psi(h^2G)-
\psi(h^2G_0)\right)\langle x\rangle^s\right\|_{L^2\to L^p}$$ $$
  \le O_N(h^2)\int_{{\bf C}_\psi}\left|{\rm Im}\,z
\right|^N\left\|
(h^2G-z)^{-1}\right\|_{L^2\to L^p}\left\|\langle x\rangle^{-s}
(h^2G_0-z)^{-1}\langle x\rangle^s
\right\|_{L^2\to L^2}L(dz).\eqno{(2.44)}$$
Now (2.34) follows from (2.36), (2.44) and the bound
$$\left\|(h^2G-z)^{-1}\right\|_{L^2\to L^p}
\le Ch^{-q}|{\rm Im}\,z|^{-\widetilde q},\eqno{(2.45)}$$
for $z\in{\bf C}_\psi$, ${\rm Im}\,z\neq 0$, with constants 
$C,\widetilde q>0$ independent of $h$ and $z$, where $q=n\left|\frac{1}{2}-
\frac{1}{p}\right|$. On the other hand, it is easy to see that
(2.45) follows from (2.41) combined with (2.38),(2.42) and (2.43).
\eproof

\section{$L^2\to L^2$ estimates for the wave group}

Given a parameter $0<h\le 1$ and a real-valued 
function $\varphi\in C_0^\infty((0,+\infty))$, denote
$$\Phi(t;h)=e^{it\sqrt{G}}\varphi(h\sqrt{G})-
e^{it\sqrt{G_0}}\varphi(h\sqrt{G_0}).$$
We will first prove the following

\begin{Theorem} Assume (1.1) fulfilled. Then we have
$$\left\|\Phi(t;h)\right\|_{L^2\to L^2}\le Ch,\quad
\forall t,\,0<h\le 1,\eqno{(3.1)}$$
with a constant $C>0$ independent of $t$ and $h$.
\end{Theorem}

{\it Proof.} We will derive (3.1) from the following

\begin{prop} Assume (1.1) fulfilled. Then, for every 
$s>1/2$ and every real-valued function 
$\varphi\in C_0^\infty((0,+\infty))$ the following estimate holds
$$\int_{-\infty}^\infty\left\|\langle x\rangle^{-s}e^{it\sqrt{G}}
\varphi(h\sqrt{G})f\right\|^2_{L^2}dt\le C\|f\|^2_{L^2},\quad\forall
 f\in L^2,\,0<h\le 1,\eqno{(3.2)}$$
with a constant $C>0$ independent of $h$ and $f$.
\end{prop}

By Duhamel's formula
$$e^{it\sqrt{G}}\varphi(h\sqrt{G})-e^{it\sqrt{G_0}}\varphi(h\sqrt{G})=
i\frac{\sin\left(t\sqrt{G_0}\right)}{\sqrt{G_0}}
\left(\sqrt{G}\varphi(h\sqrt{G})-\sqrt{G_0}\varphi(h\sqrt{G})\right)
$$ $$-\int_0^t\frac{\sin\left((t-\tau)
\sqrt{G_0}\right)}{\sqrt{G_0}}Ve^{i\tau\sqrt{G}}\varphi(h\sqrt{G})
d\tau,\eqno{(3.3)}$$
we obtain
$$\Phi(t;h)=\Phi_1(t;h)+h\Phi_2(t;h),\eqno{(3.4)}$$
where 
 $$\Phi_1(t;h)=
\left(\varphi_1(h\sqrt{G})-\varphi_1(h\sqrt{G_0})\right)
e^{it\sqrt{G}}\varphi(h\sqrt{G})$$ $$+\varphi_1(h\sqrt{G_0})
e^{it\sqrt{G_0}}\left(\varphi(h\sqrt{G})-\varphi(h\sqrt{G_0})\right)$$
$$-i\varphi_1(h\sqrt{G_0})\sin\left(t\sqrt{G_0}\right)
\left(\varphi(h\sqrt{G})-\varphi(h\sqrt{G_0})\right)$$ $$
+i\widetilde\varphi_1(h\sqrt{G_0})\sin\left(t\sqrt{G_0}\right)
\left(\widetilde\varphi(h\sqrt{G})-
\widetilde\varphi(h\sqrt{G_0})\right),$$
 $$\Phi_2(t;h)=-\int_0^t\widetilde\varphi_1(h\sqrt{G_0})\sin\left((t-\tau)
\sqrt{G_0}\right)Ve^{i\tau\sqrt{G}}\varphi(h\sqrt{G})d\tau,$$
where $\varphi_1\in C_0^\infty((0,+\infty))$ is a real-valued
function such that 
$\varphi_1\varphi\equiv \varphi$, $\widetilde\varphi(\sigma)=
\sigma\varphi(\sigma)$, $\widetilde\varphi_1(\sigma)=
\sigma^{-1}\varphi_1(\sigma)$. For all nontrivial $f,g\in L^2$, 
in view of (2.31), (3.2) and (3.4), 
we have with $0<s-1/2\ll 1$, $\forall\gamma>0$,
$$\left|\left\langle \Phi(t;h)f,g\right\rangle\right|\le 
O(h^2)\|f\|_{L^2}\|g\|_{L^2}$$ $$+
O(h)\int_{-\infty}^\infty\left|\left\langle \langle x\rangle^{s}V
 e^{i\tau\sqrt{G}}\varphi(h\sqrt{G})f,\langle x\rangle^{-s}
\sin\left((t-\tau)\sqrt{G_0}\right)\widetilde\varphi_1(h\sqrt{G_0})g
\right\rangle\right|d\tau$$ $$
\le O(h^2)\|f\|_{L^2}\|g\|_{L^2}+
O(h)\gamma\int_{-\infty}^\infty\left\|\langle x\rangle^{-s}
 e^{i\tau\sqrt{G}}\varphi(h\sqrt{G})f\right\|^2_{L^2}d\tau$$ 
 $$+O(h)\gamma^{-1}\int_{-\infty}^\infty\left\|\langle x\rangle^{-s}
\sin\left(\tau\sqrt{G_0}\right)\widetilde\varphi_1(h\sqrt{G_0})g
\right\|^2_{L^2}d\tau$$ $$\le  O(h^2)\|f\|_{L^2}\|g\|_{L^2}+
O(h)\gamma\|f\|^2_{L^2}+
 O(h)\gamma^{-1}\|g\|^2_{L^2}\le 
 O(h)\|f\|_{L^2}\|g\|_{L^2},\eqno{(3.5)}$$
if we choose $\gamma=\|g\|_{L^2}/\|f\|_{L^2}$, 
which clearly implies (3.1).
\eproof

{\it Proof of Proposition 3.2.} 
Denote by ${\cal H}$ the Hilbert space
$L^2({\bf R};L^2)$. Clearly, (3.2) is equivalent to the fact that
the operator ${\cal A}_h:L^2\to{\cal H}$ defined by
$$\left({\cal A}_hf\right)(x,t)=\langle x\rangle^{-s}e^{it\sqrt{G}}
\varphi(h\sqrt{G})f$$
is bounded uniformly in $h$. Observe that the adjoint
${\cal A}_h^*:{\cal H}\to L^2$ is defined by
$${\cal A}_h^*f=\int_{-\infty}^\infty e^{-i\tau\sqrt{G}}\varphi(h\sqrt{G})
\langle x\rangle^{-s}f(\tau,x)d\tau,$$
so we have, $\forall f,g\in{\cal H}$,
$$\left\langle {\cal A}_h{\cal A}_h^*f,g\right\rangle_{{\cal H}}=
\int_{-\infty}^\infty\left\langle \rho(t,\cdot),g(t,\cdot)
\right\rangle_{L^2}dt,\eqno{(3.6)}$$
where
$$\rho(t,x)=\int_{-\infty}^\infty \langle x\rangle^{-s}
e^{i(t-\tau)\sqrt{G}}\varphi^2(h\sqrt{G})
\langle x\rangle^{-s}f(\tau,\cdot)d\tau.$$
Hence, for the Fourier transform, $\hat\rho(\lambda,x)$, of
$\rho(t,x)$ with respect to the variable $t$ we have
$$\hat\rho(\lambda,x)=Q(\lambda)\hat f(\lambda,x),\eqno{(3.7)}$$
where $Q(\lambda)$ is the Fourier transform of the operator
$$\langle x\rangle^{-s}
e^{it\sqrt{G}}\varphi^2(h\sqrt{G})\langle x\rangle^{-s}.$$
On the other hand, the formula
$$e^{it\sqrt{G}}\varphi^2(h\sqrt{G})=\frac{1}{\pi i}\int_{-\infty}^\infty
e^{it\lambda}\varphi^2(h\lambda)\left(R^+(\lambda)-R^-(\lambda)\right)
\lambda d\lambda,\eqno{(3.8)}$$
where
$$R^\pm(\lambda)=\lim_{\varepsilon\to 0^+}(G-\lambda^2\pm i
\varepsilon)^{-1}:\langle x\rangle^{-s}L^2\to 
\langle x\rangle^{s}L^2,\quad s>1/2,$$
shows that
$$Q(\lambda)=(\pi i)^{-1}\lambda\varphi^2(h\lambda)
\langle x\rangle^{-s}\left(R^+(\lambda)-R^-(\lambda)\right)
\langle x\rangle^{-s}.\eqno{(3.9)}$$
Note that the limit exists in view of the limiting absorption principle.
Moreover, we have the following

\begin{lemma} Assume (1.1) fulfilled. Then we have
$$\left\|\langle x\rangle^{-1/2-\epsilon}R^\pm(\lambda)
\langle x\rangle^{-1/2-\epsilon}\right\|_{L^2\to L^2}\le C\lambda^{-1},
\quad \lambda\ge\lambda_0,\eqno{(3.10)}$$
for every $\lambda_0>0$, $0<\epsilon\ll 1$, with a constant
$C>0$ independent of $\lambda$. 
\end{lemma}

{\it Proof.} The estimate (3.10) is well known to 
hold for the free operator $G_0$, i.e. we have
$$\left\|\langle x\rangle^{-1/2-\epsilon}R_0^\pm(\lambda)
\langle x\rangle^{-1/2-\epsilon}\right\|_{L^2\to L^2}\le C\lambda^{-1},
\quad\forall\lambda>0.\eqno{(3.11)}$$
To show that it holds for the perturbed operator as well, 
we will take advantage of the identity
$$\langle x\rangle^{-s}R^\pm(\lambda)\langle x\rangle^{-s_1}
\left(1+K^\pm(\lambda)\right)=\langle x\rangle^{-s}
R^\pm_0(\lambda)\langle x\rangle^{-s_1},\eqno{(3.12)}$$
where
$$K^\pm(\lambda)=\langle x\rangle^{s_1}V
R^\pm_0(\lambda)\langle x\rangle^{-s}.$$
By (3.11), we have with $1/2<s_1\le\delta-1/2$,
$$\left\|K^\pm(\lambda)\right\|_{L^2\to L^2}\le C\lambda^{-1},
\quad\forall\lambda>0.\eqno{(3.13)}$$
Hence, there exists $\lambda_0>0$ so that we have
$$\left\|\left(1+K^\pm(\lambda)\right)^{-1}
\right\|_{L^2\to L^2}\le Const,
\quad\forall\lambda\ge\lambda_0.\eqno{(3.14)}$$
Moreover, since 
$G$ has no strictly positive resonances, (3.14) holds for any
$\lambda_0>0$. Now (3.10) follows from (3.11), (3.12) and (3.14). 
\eproof

By (3.9) and (3.10) we conclude
$$\|Q(\lambda)\|_{L^2\to L^2}\le C\eqno{(3.15)}$$
with a constant $C>0$ independent of $\lambda$ and $h$. By (3.7) and
(3.15),
$$\|\hat\rho(\lambda,\cdot)\|_{L^2}\le C\|\hat f(\lambda,\cdot)\|_{L^2},
\eqno{(3.16)}$$
which together with (3.6) leads to
$$\left|\left\langle {\cal A}_h{\cal A}_h^*f,g\right
\rangle_{{\cal H}}\right|=\left|
\int_{-\infty}^\infty\left\langle \hat\rho
(\lambda,\cdot),\hat g(\lambda,\cdot)
\right\rangle_{L^2}d\lambda\right|$$ $$\le C
\int_{-\infty}^\infty\|\hat f
(\lambda,\cdot)\|_{L^2}\|\hat g(\lambda,\cdot)\|_{L^2}d\lambda$$ 
 $$\le C\gamma\int_{-\infty}^\infty\|\hat f
(\lambda,\cdot)\|_{L^2}^2d\lambda+
C\gamma^{-1}\int_{-\infty}^\infty\|\hat g
(\lambda,\cdot)\|_{L^2}^2d\lambda$$ 
 $$=C\gamma\|f\|^2_{{\cal H}}+C\gamma^{-1}\|g\|^2_{{\cal H}}=2C
\|f\|_{{\cal H}}\|g\|_{{\cal H}},\eqno{(3.17)}$$
if we take $\gamma=\|g\|_{{\cal H}}/\|f\|_{{\cal H}}$, with a constant
$C>0$ independent of $h$, $f$ and $g$. It follows from (3.17) that
the operator ${\cal A}_h{\cal A}_h^*:{\cal H}\to{\cal H}$ is bounded
uniformly in $h$, and hence so is the operator ${\cal A}_h:L^2\to
{\cal H}$. This clearly proves (3.2).
\eproof

In what follows in this section we will prove the following

\begin{Theorem} Assume (1.1) fulfilled. Then, for every real-valued function 
$\varphi\in C_0^\infty((0,+\infty))$ and every 
$0\le s\le(n-1)/2$, $0<\epsilon\ll 1$, we have
$$\left\|\langle x\rangle^{-s-\epsilon}
e^{it\sqrt{G}}\varphi(h\sqrt{G})\langle x\rangle^{-s-\epsilon}
\right\|_{L^2\to L^2}\le C\langle t\rangle^{-s},\quad 
\forall t,\,0<h\le 1,\eqno{(3.18)}$$
with a constant $C>0$ independent of $t$ and $h$.
\end{Theorem}

{\it Proof.} We will derive (3.18) from the following estimates 

\begin{prop} Assume (1.1) fulfilled. Then, for every real-valued function 
$\varphi\in C_0^\infty((0,+\infty))$ and every 
$0\le s\le(n-1)/2$, $0<\epsilon\ll 1$, $0<h\le 1$, we have
$$\left\|\langle x\rangle^{-1/2-s-\epsilon}
e^{it\sqrt{G}}\varphi(\sqrt{G})\langle x\rangle^{-1/2-s-\epsilon}
\right\|_{L^2\to L^2}\le C\langle t\rangle^{-s},\quad 
\forall t,\eqno{(3.19)}$$
$$\int_{-\infty}^\infty|t|^{2s}\left\|\langle x\rangle^{-1/2-s-\epsilon}
e^{it\sqrt{G}}
\varphi(h\sqrt{G})\langle x\rangle^{-1/2-s-\epsilon}
f\right\|^2_{L^2}dt\le C\|f\|^2_{L^2},\quad\forall f\in L^2,\eqno{(3.20)}$$
with a constant $C>0$ independent of $t$, $h$ and $f$. 
\end{prop}

We will first show that (3.18) with $0\le s\le(n-1)/2$ follows from
(3.18) with $s=(n-1)/2$. To this end, consider the operator-valued function
$$P(z)=\langle t\rangle^{z}\langle x\rangle^{-z-\epsilon}
e^{it\sqrt{G}}\varphi(h\sqrt{G})\langle x\rangle^{-z-\epsilon},\quad
 z\in{\bf C}.$$ 
Clearly, $P(z)$ is analytic in $z$ for ${\rm Re}\,z\ge 0$ with values in
${\cal L}(L^2)$ and satisfies the trivial bounds
$$\left\|P(z)\right\|_{L^2\to L^2}\le C\langle t\rangle^{{\rm Re}\,z},
\quad 0\le {\rm Re}\,z\le(n-1)/2.\eqno{(3.21)}$$
On the other hand, supposing that (3.18) holds with $s=(n-1)/2$ implies
$$\left\|P(z)\right\|_{L^2\to L^2}\le C\eqno{(3.22)}$$
on ${\rm Re}\,z=0$ and ${\rm Re}\,z=(n-1)/2$, with a constant $C>0$
independent of $z$, $t$ and $h$. It follows now from (3.21), (3.22)
and the Phragm\`en-Lindel\"of principle that (3.22) holds for
$0\le {\rm Re}\,z\le(n-1)/2$, which is the desired result.

Using (1.1), (2.1) and (2.28), we obtain 
$$\left\|\langle x\rangle^{-s-\epsilon}
\Phi_1(t;h)\langle x\rangle^{-s_1-\epsilon}f
\right\|_{L^2}$$ $$\le O(h^2)
\left\|\langle x\rangle^{-(n+1)/2-\epsilon}
\Phi(t;h)\langle x\rangle^{-s_1-\epsilon}f
\right\|_{L^2}+O(h^2)\langle t\rangle^{-s}\|f\|_{L^2},
\quad\forall f\in L^2,\eqno{(3.23)}$$
for all $s_1\ge s\ge 0$. Using (1.1), (2.1) and (3.20), we have
$\forall f,g\in L^2$, with $s=(n-1)/2$,
$$\langle t\rangle^s\left|\left\langle \langle x\rangle^{-s-\epsilon}
\Phi_2(t;h)\langle x\rangle^{-s-1/2-\epsilon}f,g
\right\rangle\right|$$ $$\le \langle t\rangle^s 
  \int_0^t\left|\left\langle \langle x\rangle^{-s-\epsilon}
\widetilde\varphi_1(h\sqrt{G_0})\sin\left((t-\tau)
\sqrt{G_0}\right)Ve^{i\tau\sqrt{G}}\varphi(h\sqrt{G}) 
\langle x\rangle^{-1/2-s-\epsilon}f,g\right\rangle\right|d\tau$$
 $$\le 
  C\int_0^{t/2}\langle t-\tau\rangle^s\left\|\langle x
\rangle^{-(n+\epsilon)/2}e^{i(t-\tau)\sqrt{G}}\varphi(h\sqrt{G}) 
\langle x\rangle^{-1/2-s-\epsilon}f\right\|_{L^2}$$
 $$\times \left\|\langle x
\rangle^{-(1+\epsilon)/2}\sin\left(\tau\sqrt{G_0}\right)
\widetilde\varphi_1(h\sqrt{G_0})\langle x\rangle^{-s-\epsilon}
g\right\|_{L^2}d\tau$$
 $$+C\int_0^{t/2}\left\|\langle x
\rangle^{-1-\epsilon}e^{i\tau\sqrt{G}}\varphi(h\sqrt{G}) 
\langle x\rangle^{-1/2-s-\epsilon}f\right\|_{L^2}$$
 $$\times\langle t-\tau\rangle^s\left\|\langle x
\rangle^{-(n-1+\epsilon)/2}\sin\left((t-\tau)\sqrt{G_0}\right)
\widetilde\varphi_1(h\sqrt{G_0})\langle x\rangle^{-s-\epsilon}
g\right\|_{L^2}d\tau$$
 $$\le C\gamma\int_{-\infty}^{\infty}\langle \tau\rangle^{2s}
\left\|\langle x
\rangle^{-1/2-s-\epsilon}e^{i\tau\sqrt{G}}\varphi(h\sqrt{G}) 
\langle x\rangle^{-1/2-s-\epsilon}f\right\|_{L^2}^2d\tau$$
 $$+C\gamma^{-1}\int_{-\infty}^{\infty}\left\|\langle x
\rangle^{-(1+\epsilon)/2}\sin\left(\tau\sqrt{G_0}\right)
\widetilde\varphi_1(h\sqrt{G_0})\langle x\rangle^{-s-\epsilon}
g\right\|_{L^2}^2d\tau$$
 $$+ C\gamma\int_{-\infty}^{\infty}
\langle\tau\rangle^{1+\epsilon}\left\|\langle x\rangle^{-1-\epsilon}
e^{i(t-\tau)\sqrt{G}}\varphi(h\sqrt{G}) 
\langle x\rangle^{-1/2-s-\epsilon}f\right\|_{L^2}^2d\tau$$
 $$+C\gamma^{-1}\int_{-\infty}^{\infty}\langle\tau\rangle^{-1-\epsilon}
\langle t-\tau\rangle^{2s}\left\|\langle x
\rangle^{-s-\epsilon}\sin\left((t-\tau)\sqrt{G_0}\right)
\widetilde\varphi_1(h\sqrt{G_0})\langle x\rangle^{-s-\epsilon}
g\right\|_{L^2}^2d\tau$$
 $$\le C\gamma\|f\|_{L^2}^2+C\gamma^{-1}\|g\|_{L^2}^2
=2C\|f\|_{L^2}\|g\|_{L^2},\eqno{(3.24)}$$
if we take $\gamma=\|g\|_{L^2}/\|f\|_{L^2}$. 
By (3.4), (3.23) and (3.24), we obtain (with $s=(n-1)/2$)
$$\left\|\langle x\rangle^{-s-\epsilon}
\Phi(t;h)\langle x\rangle^{-s-1/2-\epsilon}f
\right\|_{L^2}$$ $$\le O(h^2)
\left\|\langle x\rangle^{-(n+1)/2-\epsilon}
\Phi(t;h)\langle x\rangle^{-s-1/2-\epsilon}f
\right\|_{L^2}+O(h)\langle t\rangle^{-s}\|f\|_{L^2}.\eqno{(3.25)}$$
Hence, there exists a constant $0<h_0<1$ so that for $0<h\le h_0$
we can absorbe the first term in the RHS of (3.25), thus obtaining the 
estimate ( for  $0<h\le h_0$)
$$\left\|\langle x\rangle^{-s-\epsilon}
\Phi(t;h)\langle x\rangle^{-1/2-s-\epsilon}
\right\|_{L^2\to L^2}\le Ch\langle t\rangle^{-s}.\eqno{(3.26)}$$
Let now $h_0\le h\le 1$. Without loss of generality we may suppose
 $h=1$. Then, by (2.1) and (3.19), the norm in the first term in
the RHS of (3.25) is upper bounded by $C\langle t\rangle^{-s}\|f\|_{L^2}$,
which again implies (3.26). By (2.1) and (3.26), we have
$$\left\|\langle x\rangle^{-s-\epsilon}
e^{it\sqrt{G}}\varphi(h\sqrt{G})\langle x\rangle^{-1/2-s-\epsilon}
\right\|_{L^2\to L^2}\le C\langle t\rangle^{-s},
\quad\forall t,\,0<h\le 1,\eqno{(3.27)}$$
where $s=(n-1)/2$. By the same interpolation argument as above, we conclude
that (3.27) holds for all $0\le s\le (n-1)/2$. We will show now that this
implies (3.18) with $s=(n-1)/2$. 
To this end, we will proceed in a way similar to that one
in Section 3 of \cite{kn:V2} which is based on the commutator relation
$$ -2\Delta+[r\partial_r,\Delta]=0,\eqno{(3.28)}$$
where $r=|x|$ is the radial variable. In the same way, one can show that
(3.27) implies
$$\left\|\langle x\rangle^{-s-\epsilon}{\cal D}_r
e^{it\sqrt{G}}\varphi(h\sqrt{G})\langle x\rangle^{-1/2-s-\epsilon}
\right\|_{L^2\to L^2}\le C\langle t\rangle^{-s},
\quad\forall t,\,0<h\le 1,\eqno{(3.29)}$$
for all $0\le s\le (n-1)/2$, where ${\cal D}_r:=\langle r\rangle^{-1}
rh\partial_r$. Furthermore, by Duhamel's formula and (3.28), 
we have the identity
$$\varphi_1(h\sqrt{G})[r\partial_r,e^{it\sqrt{G}}]\varphi(h\sqrt{G}) 
 =i\varphi_1(h\sqrt{G})\frac{\sin\left(t\sqrt{G}\right)}
{\sqrt{G}}[r\partial_r,\sqrt{G}]\varphi(h\sqrt{G})$$ 
$$-\int_0^t\varphi_1(h\sqrt{G})\frac{\sin\left((t-\tau)\sqrt{G}\right)}
{\sqrt{G}}[r\partial_r,G]e^{i\tau\sqrt{G}}\varphi(h\sqrt{G})d\tau$$
 $$=i\varphi_1(h\sqrt{G})\frac{\sin\left(t\sqrt{G}\right)}
{\sqrt{G}}[r\partial_r,\sqrt{G}]\varphi(h\sqrt{G})$$ 
 $$-\int_0^t\varphi_1(h\sqrt{G})\frac{\sin\left((t-\tau)\sqrt{G}\right)}
{\sqrt{G}}\left(2G-2V+r\partial_rV-Vr\partial_r\right)
e^{i\tau\sqrt{G}}\varphi(h\sqrt{G})d\tau$$
 $$=i\varphi_1(h\sqrt{G})\frac{\sin\left(t\sqrt{G}\right)}
{\sqrt{G}}r\partial_r\sqrt{G}\varphi(h\sqrt{G})
-i\varphi_1(h\sqrt{G})\sin\left(t\sqrt{G}\right)
r\partial_r\varphi(h\sqrt{G})$$ 
 $$+ite^{it\sqrt{G}}\sqrt{G}\varphi(h\sqrt{G})+2i\sin\left(t\sqrt{G}
\right)\varphi(h\sqrt{G})$$ 
 $$-\int_0^t\varphi_1(h\sqrt{G})\frac{\sin\left((t-\tau)\sqrt{G}\right)}
{\sqrt{G}}\left(-V+\partial_rrV-Vr\partial_r\right)
e^{i\tau\sqrt{G}}\varphi(h\sqrt{G})d\tau,$$
where the functions $\varphi$ and $\varphi_1$ are as above. Set 
$\varphi_2(\sigma)=\sigma^{-1}\varphi(\sigma)$, $\widetilde\varphi_1
(\sigma)=\sigma^{-1}\varphi_1(\sigma)$. From the above identity we get
(with $s=(n-1)/2$)
$$it\langle x\rangle^{-s-\epsilon}
e^{it\sqrt{G}}\varphi(h\sqrt{G})\langle x\rangle^{-s-\epsilon}=-2ih
\langle x\rangle^{-s-\epsilon}
\sin\left(t\sqrt{G}\right)\varphi_2(h\sqrt{G})\langle x\rangle^{-s-\epsilon}$$
 $$+i\langle x\rangle^{-s-\epsilon}\widetilde\varphi_1(h\sqrt{G})
\sin\left(t\sqrt{G}\right)({\cal D}_r^*+O(h))
\langle x\rangle^{-s+1-\epsilon}\left(
\langle x\rangle^{s+\epsilon}
\varphi(h\sqrt{G})\langle x\rangle^{-s-\epsilon}\right)$$
 $$-i\langle x\rangle^{-s-\epsilon}\varphi_1(h\sqrt{G})
\sin\left(t\sqrt{G}\right)({\cal D}_r^*+O(h))
\langle x\rangle^{-s+1-\epsilon}\left(
\langle x\rangle^{s+\epsilon}
\varphi_2(h\sqrt{G})\langle x\rangle^{-s-\epsilon}\right)$$
 $$+\left(\langle x\rangle^{-s-\epsilon}\varphi_1(h\sqrt{G})
\langle x\rangle^{s+\epsilon}\right)\langle x\rangle^{-s+1-\epsilon}
{\cal D}_re^{it\sqrt{G}}\varphi_2(h\sqrt{G})\langle x\rangle^{-s-\epsilon}$$
 $$+\langle x\rangle^{-s-\epsilon}\varphi_1(h\sqrt{G})
e^{it\sqrt{G}}({\cal D}_r^*+O(h))\langle x\rangle^{-s+1-\epsilon}\left(
\langle x\rangle^{s+\epsilon}
\varphi_2(h\sqrt{G})\langle x\rangle^{-s-\epsilon}\right)$$
 $$+h^2\int_0^t\langle x\rangle^{-s-\epsilon}
\widetilde\varphi_1(h\sqrt{G})\sin\left((t-\tau)\sqrt{G}\right)
Ve^{i\tau\sqrt{G}}\varphi_2(h\sqrt{G})\langle x\rangle^{-s-\epsilon}d\tau$$
 $$+h\int_0^t\langle x\rangle^{-s-\epsilon}
\widetilde\varphi_1(h\sqrt{G})\sin\left((t-\tau)\sqrt{G}\right)
({\cal D}_r^*+O(h))\langle r\rangle
Ve^{i\tau\sqrt{G}}\varphi_2(h\sqrt{G})\langle x\rangle^{-s-\epsilon}d\tau$$
 $$+h\int_0^t\langle x\rangle^{-s-\epsilon}
\widetilde\varphi_1(h\sqrt{G})\sin\left((t-\tau)\sqrt{G}\right)
V\langle r\rangle {\cal D}_r
e^{i\tau\sqrt{G}}\varphi_2(h\sqrt{G})\langle x\rangle^{-s-\epsilon}d\tau.
\eqno{(3.30)}$$
By (2.27), (3.27) and (3.29), we have that the $L^2\to L^2$ norm of each
of the first five terms in the RHS of (3.30) is upper bounded by
$O(\langle t\rangle^{-s+1})$. In the same way, taking into acount (1.1),
one can easily see that the $L^2\to L^2$ norm of each integral in
the RHS of (3.30) is upper  bounded by $O(\langle t\rangle^{-s+1/2})$.
Thus, (3.30) implies (3.18) with $s=(n-1)/2$, and hence with all
$0\le s\le (n-1)/2$. 
\eproof

{\it Proof of Proposition 3.5.} We will derive (3.19) from the following

\begin{lemma} Assume (1.1) fulfilled and let $0\le s\le (n-1)/2$.  
Let also $m\ge 0$ denote the bigest integer $\le s$ and set
$\mu=s-m$. Then, the
operator-valued function
$${\cal R}_s^\pm(\lambda)=\lambda \langle x\rangle^{-1/2-s-\epsilon}
R^\pm(\lambda)\langle x\rangle^{-1/2-s-\epsilon}:L^2\to L^2$$
is $C^m$ in $\lambda$ for $\lambda>0$, 
$\partial_\lambda^m{\cal R}_s^\pm$ is H\"older of order $\mu$, and
satisfies the estimates
$$\left\|\partial_\lambda^j{\cal R}_s^\pm(\lambda)
\right\|_{L^2\to L^2}\le C,\quad \lambda\ge\lambda_0,\,0\le j\le m,
\eqno{(3.31)}$$
$$\left\|\partial_\lambda^m{\cal R}_s^\pm(\lambda_2)-
\partial_\lambda^m{\cal R}_s^\pm(\lambda_1)
\right\|_{L^2\to L^2}\le C|\lambda_2-\lambda_1|^\mu,\quad 
\lambda_2>\lambda_1\ge\lambda_0,\eqno{(3.32)}$$
for every $\lambda_0>0$, 
with a constant $C>0$ independent of $\lambda$, 
$\lambda_1$ and $\lambda_2$. 
\end{lemma}

{\it Proof.} We will derive (3.31) and (3.32) from the fact that 
for every integer $j\ge 0$ the operator-valued function 
$${\cal R}_{0,s}^\pm(\lambda)=\lambda \langle x\rangle^{-1/2-s-\epsilon}
R^\pm_0(\lambda)\langle x\rangle^{-1/2-s-\epsilon}:L^2\to L^2$$
is $C^j$ in $\lambda$ if $s\ge j$, and satisfies the bound
$$\left\|\partial_\lambda^j{\cal R}_{0,j}^\pm(\lambda)
\right\|_{L^2\to L^2}\le C_j,\quad \forall\lambda>0.\eqno{(3.33)}$$
Let $s$, $m$ and $\mu$ be as in Lemma 3.6. We will show that (3.33) implies
$$\left\|\partial_\lambda^m{\cal R}_{0,s}^\pm(\lambda_2)-
\partial_\lambda^m{\cal R}_{0,s}^\pm(\lambda_1)
\right\|_{L^2\to L^2}\le C|\lambda_2-\lambda_1|^\mu,\quad 
\lambda_2>\lambda_1>0.\eqno{(3.34)}$$
To this end, fix $\lambda_2>\lambda_1$ and consider the operator-valued
function
$${\cal P}_\pm(z)=|\lambda_2-\lambda_1|^{-z}\left(
\partial_\lambda^m{\cal R}_{0,z}^\pm(\lambda_2)-
\partial_\lambda^m{\cal R}_{0,z}^\pm(\lambda_1)\right),\quad
 z\in{\bf C}.$$
In view of (3.33), ${\cal P}_\pm(z)$ is analytic in $z$ for ${\rm Re}\,z\ge 0$
with values in ${\cal L}(L^2)$ and satisfies the bounds
$$\|{\cal P}_\pm(z)\|_{L^2\to L^2}\le C_1,
\quad 0\le {\rm Re}\,z\le 1,\eqno{(3.35)}$$
with a constant $C_1>0$ independent of $z$ but 
depending on $\lambda_1$ and $\lambda_2$,
while on ${\rm Re}\,z=0$ and ${\rm Re}\,z=1$ we have a better bound
$$\|{\cal P}_\pm(z)\|_{L^2\to L^2}\le C,\eqno{(3.36)}$$
with a constant $C>0$ independent of $z$, $\lambda_1$ and $\lambda_2$. 
By (3.35), (3.36) and the Phragm\`en-Lindel\"of principle 
we conclude that (3.36) holds for $0\le {\rm Re}\,z\le 1$. 
In particular, it holds for $z=\mu$, which proves (3.34).

To prove (3.31) we differentiate $j$ times the identity (3.12)
to obtain
$$\frac{d^j{\cal R}_j^\pm}{d\lambda^j}(\lambda)
\left(1+K_j^\pm(\lambda)\right)=
\frac{d^j{\cal R}_{0,j}^\pm}{d\lambda^j}(\lambda)
 -\sum_{k=0}^j\beta_k\langle x\rangle^{-j+k}
\frac{d^k{\cal R}_k^\pm}{d\lambda^k}(\lambda)
\langle x\rangle^{-j+k}
\frac{d^{j-k}K_j^\pm}{d\lambda^{j-k}}(\lambda),
\eqno{(3.37)}$$
where 
$$K_j^\pm(\lambda)=\langle x\rangle^{1/2+j+\epsilon}VR_0^\pm(\lambda)
\langle x\rangle^{-1/2-j-\epsilon},\quad 0\le j\le (n-1)/2.$$
As in the proof of Lemma 3.3 above, we conclude that (3.14) still
holds for $K_j^\pm(\lambda)$. Therefore, (3.31) follows from
(3.37) combined with (3.10) by induction in $j$.

To prove (3.32) observe that (3.37) leads to the identity
$$\left(\frac{d^m{\cal R}_s^\pm}{d\lambda^m}(\lambda_2)-
\frac{d^m{\cal R}_s^\pm}{d\lambda^m}(\lambda_1)\right)
\left(1+K_m^\pm(\lambda_2)\right)=
\frac{d^m{\cal R}_{0,s}^\pm}{d\lambda^m}(\lambda_2)-
\frac{d^m{\cal R}_{0,s}^\pm}{d\lambda^m}(\lambda_1)$$ 
 $$-\sum_{k=0}^m\beta_k\langle x\rangle^{-m+k}
\frac{d^k{\cal R}_s^\pm}{d\lambda^k}(\lambda_1)
\left(\langle x\rangle^{-m+k}
\frac{d^{m-k}K_m^\pm}{d\lambda^{m-k}}(\lambda_2)
-\langle x\rangle^{-m+k}\frac{d^{m-k}K_m^\pm}{d\lambda^{m-k}}
(\lambda_1)\right)$$ $$
-\sum_{k=0}^{m-1}\beta_k\langle x\rangle^{-m+k}\left(
\frac{d^k{\cal R}_s^\pm}{d\lambda^k}(\lambda_2)
-\frac{d^k{\cal R}_s^\pm}{d\lambda^k}(\lambda_1)\right)
\langle x\rangle^{-m+k}
\frac{d^{m-1-k}K_m^\pm}{d\lambda^{m-1-k}}(\lambda_2).
\eqno{(3.38)}$$
As above, one has that (3.14) still holds for $K_m^\pm(\lambda_2)$,
so one can easily derive (3.32) from (3.31), (3.33), (3.34) and (3.38).
\eproof

We are going to use (3.8) with 
$h=1$ and $\varphi^2$ replaced by $\varphi$. Set 
$$T(\lambda)=T^+(\lambda)-T^-(\lambda),\quad 
T^\pm(\lambda)=(\pi i)^{-1}
{\cal R}_s^\pm(\lambda),\eqno{(3.39)}$$ 
and choose a real-valued
function $\phi\in C_0^\infty([1/3,1/2])$, $\phi\ge 0$, such that 
$\int\phi(\sigma)d\sigma=1$. Then the function
$$T^\pm_\theta(\lambda)=\theta^{-1}\int T^\pm(\lambda+\sigma)\phi(\sigma/
\theta)d\sigma,\quad 0<\theta\le1,$$
is $C^\infty$ in $\lambda$ with values in ${\cal L}(L^2)$ and, 
in view of (3.31) and (3.32), satisfies the estimates 
$$\|\partial_\lambda^jT^\pm_\theta(\lambda)\|_{L^2\to L^2}\le 
C,\quad 0\le j\le m,\eqno{(3.40)}$$
 $$\|\partial_\lambda^mT^\pm_\theta(\lambda)-\partial_\lambda^m
T^\pm(\lambda)\|_{L^2\to L^2}\le \theta^{-1}
\int\|\partial_\lambda^mT^\pm(\lambda+\sigma)-\partial_\lambda^m
T^\pm(\lambda)\|_{L^2\to L^2}\phi(\sigma/\theta)d\sigma$$
 $$\le C\theta^{-1}\int\sigma^{\mu}\phi(\sigma/\theta)d\sigma\le 
 C\theta^{\mu},\eqno{(3.41)}$$ 
$$\|\partial_\lambda^jT^\pm_\theta(\lambda)-\partial_\lambda^j
T^\pm(\lambda)\|_{L^2\to L^2}\le C\theta,\quad 0\le j\le m-1,
\eqno{(3.42)}$$
 $$\left\|\partial_\lambda^{m+1}T^\pm_\theta(\lambda)
\right\|_{L^2\to L^2}\le 
\theta^{-2}\int\|\partial_\lambda^mT^\pm(\lambda+\sigma)-
\partial_\lambda^mT^\pm(\lambda)\|_{L^2\to L^2}
|\phi'(\sigma/\theta)|d\sigma$$ $$
\le C\theta^{-2}\int\sigma^{\mu}|\phi'(\sigma/\theta)|
d\sigma\le C\theta^{-1+\mu}.\eqno{(3.43)}$$ 
Integrating by parts $m$ times and using 
(3.41) and (3.42), we get
$$\left\|\int_0^\infty
 e^{it\lambda}\varphi(\lambda)\left(T^\pm_\theta(\lambda)-
T^\pm(\lambda)\right) d\lambda\right\|_{L^2\to L^2}$$ 
 $$=\left\|t^{-m}\int_0^\infty
 e^{it\lambda}\frac{d^m}{d\lambda^m}\left(
\varphi(\lambda)\left(T^\pm_\theta(\lambda)-
T^\pm(\lambda)\right)\right) d\lambda\right\|_{L^2\to L^2}
\le C\theta^{\mu}|t|^{-m}.\eqno{(3.44)}$$
Similarly, integrating by parts $m+1$ times and using (3.40) and 
(3.43), we get
$$\left\|\int_0^\infty e^{it\lambda}\varphi(\lambda)
T^\pm_\theta(\lambda) d\lambda\right\|_{L^2\to L^2}$$
 $$=\left\|t^{-m-1}\int_0^\infty
 e^{it\lambda}\frac{d^{m+1}}{d\lambda^{m+1}}\left(\varphi(\lambda)
T^\pm_\theta(\lambda)\right)d\lambda\right\|_{L^2\to L^2}
\le C\theta^{-1+\mu}|t|^{-m-1}.\eqno{(3.45)}$$
By (3.44) and (3.45), 
$$\left\|\int_0^\infty e^{it\lambda}\varphi(\lambda)
T(\lambda) d\lambda\right\|_{L^2\to L^2}\le 
 C\theta^{\mu}|t|^{-m}\left(1+|t|^{-1}\theta^{-1}
\right)\le C|t|^{-m-\mu},\eqno{(3.46)}$$
if we take $\theta=|t|^{-1}$, which clearly implies (3.19).\\

In what follows in this section we will derive (3.20) 
from Lemma 3.6. Let $0\le s\le (n-1)/2$ and let $m\ge 0$ be the
bigest integer $\le s$. Remark that the function 
$\partial_\lambda^m{\cal R}_s^\pm$ satisfies (3.32) with
$\mu=s-m+\epsilon/2$. Consequently, the estimates (3.41)
and (3.43) are valid with $\mu=s-m+\epsilon/2$. Let 
$\phi_+\in C^\infty({\bf R})$, $\phi_+(t)=0$ for $t\le 1$,
$\phi_+(t)=1$ for $t\ge 2$. Given any function $f\in L^2$, set
$$u(t;h)=\langle x\rangle^{-1/2-s-\epsilon}\phi_+(t)e^{it\sqrt{G}}
\varphi(h\sqrt{G})\langle x\rangle^{-1/2-s-\epsilon}f.$$
We have
$$\left(\partial_t^2+G\right)\langle x\rangle^{1/2+s+\epsilon}u(t;h)=
2h^{-1}\phi'_+(t)e^{it\sqrt{G}}
\widetilde\varphi(h\sqrt{G})\langle x\rangle^{-1/2-s-\epsilon}f$$ $$
+\phi''_+(t)e^{it\sqrt{G}}
\varphi(h\sqrt{G})\langle x\rangle^{-1/2-s-\epsilon}f=:h^{-1}
\langle x\rangle^{-1/2-s-\epsilon}v(t;h),\eqno{(3.47)}$$
where $\widetilde\varphi(\sigma)=\sigma\varphi(\sigma)$. Clearly,
the support of the function $v(t;h)$ with respect to the variable $t$
is contained in the interval $[1,2]$, and by (2.27) we have 
$$\|v(t;h)\|_{L^2}\le C\|f\|_{L^2},\quad 1\le t\le 2,\eqno{(3.48)}$$
with a constant $C>0$ independent of $t$, $h$ and $f$. Using
Duhamel's formula we deduce from (3.47), for $t>0$, 
$$u(t;h)=\int_0^t\langle x\rangle^{-1/2-s-\epsilon}
\widetilde\varphi_1(h\sqrt{G})\sin\left((t-\tau)
\sqrt{G}\right)\langle x\rangle^{-1/2-s-\epsilon}
v(\tau;h)d\tau,\eqno{(3.49)}$$
where the functions $\varphi_1$ and $\widetilde\varphi_1$ are as above. 
It follows from (3.49) that the Fourier transforms of the functions $u(t;h)$
and $v(t;h)$ satisfy the identity
$$\hat u(\lambda;h)=Q^+(\lambda)\hat v(\lambda;h),
\eqno{(3.50)}$$
where $Q^+(\lambda)$ is the Fourier transform of the operator
$$\langle x\rangle^{-1/2-s-\epsilon}
\widetilde\varphi_1(h\sqrt{G})\eta_+(t)\sin\left(t\sqrt{G}\right)
\langle x\rangle^{-1/2-s-\epsilon},$$
$\eta_+$ being the characteristic function of the interval $[0,+\infty)$. 
It is easy to see that
$$Q^+(\lambda)={\cal B}(h)T^+(\lambda),
\eqno{(3.51)}$$
where the operator
$${\cal B}(h)=\langle x\rangle^{-1/2-s-\epsilon}
\widetilde\varphi_1(h\sqrt{G})\langle x\rangle^{1/2+s+\epsilon}:L^2\to L^2$$
is bounded uniformly in $h$ in view of (2.27). 
Set $Q^+_\theta(\lambda)={\cal B}(h)T^+_\theta(\lambda)$ 
and define the function $u_\theta(t;h)$ via the formula
$$\hat u_\theta(\lambda;h)=Q^+_\theta(\lambda)\hat v(\lambda;h).$$
Using (3.40)-(3.43) (with $\mu=s-m+\epsilon/2$) together with the 
Plancherel identity and (3.48), we obtain
$$\int_{-\infty}^\infty|t|^{2m}\|u_\theta(t;h)-u(t;h)\|_{L^2}^2dt
=\int_{-\infty}^\infty\|\partial_\lambda^m(\hat u_\theta(\lambda;h)
-\hat u(\lambda;h))\|_{L^2}^2d\lambda$$ 
 $$\le C\sum_{k=0}^m\int_{-\infty}^\infty\left\|\left(
\partial_\lambda^kQ^+_\theta(\lambda)-\partial_\lambda^k
Q^+(\lambda)\right)
\partial_\lambda^{m-k}\hat v(\lambda;h)\right\|_{L^2}^2d\lambda$$ 
 $$\le C\theta^{2\mu}\sum_{k=0}^m\int_{-\infty}^\infty\left\|
\partial_\lambda^{m-k}\hat v(\lambda;h)\right\|_{L^2}^2d\lambda
=C\theta^{2\mu}\sum_{k=0}^m\int_{-\infty}^\infty
|t|^{2m-2k}\|v(t;h)\|_{L^2}^2dt$$
 $$\le C\theta^{2\mu}\int_1^2\|v(t;h)\|_{L^2}^2dt
\le  C\theta^{2\mu}\|f\|_{L^2}^2,\eqno{(3.52)}$$
with a constant $C>0$ independent of $h$, $\theta$ and $f$.
By (3.52) we get, $\forall A\ge 1$,
$$\int_A^{2A}\|u_\theta(t;h)-u(t;h)\|_{L^2}^2dt
\le  CA^{-2m}\theta^{2\mu}\|f\|_{L^2}^2.\eqno{(3.53)}$$
In the same way, we obtain
$$\int_{-\infty}^\infty|t|^{2m+2}\|u_\theta(t;h)\|_{L^2}^2dt
=\int_{-\infty}^\infty\|\partial_\lambda^{m+1}\hat u_\theta(\lambda;h)
\|_{L^2}^2d\lambda$$ 
 $$\le C\sum_{k=0}^{m+1}\int_{-\infty}^\infty\left\|
\partial_\lambda^kQ^+_\theta(\lambda)
\partial_\lambda^{m+1-k}\hat v(\lambda;h)\right\|_{L^2}^2d\lambda$$ 
 $$\le C\theta^{-2+2\mu}\sum_{k=0}^{m+1}\int_{-\infty}^\infty\left\|
\partial_\lambda^{m+1-k}\hat v(\lambda;h)\right\|_{L^2}^2d\lambda$$ $$
=C\theta^{-2+2\mu}\sum_{k=0}^{m+1}\int_{-\infty}^\infty
|t|^{2m+2-2k}\|v(t;h)\|_{L^2}^2dt$$
 $$\le C\theta^{-2+2\mu}\int_1^2\|v(t;h)\|_{L^2}^2dt
\le  C\theta^{-2+2\mu}\|f\|_{L^2}^2,\eqno{(3.54)}$$
with a constant $C>0$ independent of $h$, $\theta$ and $f$.
By (3.54) we get, $\forall A\ge 1$,
$$\int_A^{2A}\|u_\theta(t;h)\|_{L^2}^2dt
\le  CA^{-2m-2}\theta^{-2+2\mu}\|f\|_{L^2}^2.\eqno{(3.55)}$$
Combining (3.53) and (3.55) leads to
$$\int_A^{2A}|t|^{2s}\|u(t;h)\|_{L^2}^2dt
\le  CA^{2\mu-\epsilon}\theta^{2\mu}
\left(1+A^{-2}\theta^{-2}\right)\|f\|_{L^2}^2\le
  CA^{-\epsilon}\|f\|_{L^2}^2,\eqno{(3.56)}$$
if we choose $\theta=A^{-1}$, where $s=m+\mu-\epsilon/2$.
By (3.56), for every integer $k\ge 0$ we have
$$\int_{2^k}^{2^{k+1}}|t|^{2s}\|u(t;h)\|_{L^2}^2dt
\le C2^{-\epsilon k}\|f\|_{L^2}^2.\eqno{(3.57)}$$
Summing up (3.57) leads to
$$\int_1^{\infty}|t|^{2s}\|u(t;h)\|_{L^2}^2dt
\le C\|f\|_{L^2}^2.\eqno{(3.58)}$$
It is easy to see now that this implies (3.20).
\eproof

\section{Proof of Theorem 1.1}

It is easy to see that Theorem 1.1 follows from the following

\begin{Theorem}
Assume (1.1) fulfilled. Then, for every $a>0$, $2\le p\le +\infty$,  
we have the estimates 
$$\left\|\Phi(t;h)\right\|_{L^{p'}\to L^p}
\le Ch^{1-\alpha n}|t|^{-\alpha (n-1)/2},\quad \forall t\neq 0,
\,0<h\le 1,\eqno{(4.1)}$$ 
where $1/p+1/p'=1$, $\alpha=1-2/p$, and, $\forall 0<\epsilon\ll 1$,
$$\left\|\Phi(t;h)\langle x\rangle^{-\alpha(n/2+\epsilon)}
\right\|_{L^2\to L^p}\le C_\epsilon h^{1-\alpha n/2}|t|^{-\alpha(n-1)/2},
\quad \forall t\neq 0,\,0<h\le 1,\eqno{(4.2)}$$
with constants $C,C_\epsilon>0$ independent of $t$ and $h$.
\end{Theorem}

Indeed, using the identity
$$\sigma^{-\alpha(n+1)/2}\chi_a(\sigma)=\int_0^1\varphi(\theta\sigma)
\theta^{\alpha(n+1)/2-1}d\theta,$$
where $\varphi(\sigma)=\sigma^{1-\alpha(n+1)/2}\chi'_a(\sigma)\in
C_0^\infty((0,+\infty))$,  we get 
$$\left\|e^{it\sqrt{G}}(\sqrt{G})^{-\alpha(n+1)/2}
\chi_a(\sqrt{G})-e^{it\sqrt{G_0}}(\sqrt{G_0})^{-\alpha(n+1)/2}
\chi_a(\sqrt{G_0})\right\|_{L^{p'}\to L^p}$$ $$\le \int_0^1\left\|
\Phi(t;\theta)\right\|_{L^{p'}\to L^p}\theta^{\alpha(n+1)/2-1}d\theta 
\le C|t|^{-\alpha(n-1)/2}\int_0^1\theta^{-\alpha(n-1)/2}d\theta
\le C|t|^{-\alpha(n-1)/2},\eqno{(4.3)}$$
as long as $\alpha(n-1)/2<1$, that is, for $2\le p<\frac{2(n-1)}{n-3}$. 
In the same way, we have 
$$\left\|e^{it\sqrt{G}}(\sqrt{G})^{-\alpha(n+1)/2}
\chi_a(\sqrt{G})\langle x\rangle^{-\alpha(n/2+\epsilon)}
-e^{it\sqrt{G_0}}(\sqrt{G_0})^{-\alpha(n+1)/2}
\chi_a(\sqrt{G_0})\langle x\rangle^{-\alpha(n/2+\epsilon)}
\right\|_{L^2\to L^p}$$ $$\le \int_0^1\left\|
\Phi(t;\theta)\langle x\rangle^{-\alpha(n/2+\epsilon)}
\right\|_{L^2\to L^p}\theta^{\alpha(n+1)/2-1}d\theta 
\le C|t|^{-\alpha(n-1)/2}\int_0^1\theta^{\alpha/2}d\theta
\le C|t|^{-\alpha(n-1)/2},\eqno{(4.4)}$$
which implies (1.4). Furthermore, taking 
$\varphi(\sigma)=\sigma^{1-\alpha(n-1)}\chi'_a(\sigma)$, 
since $0\le \alpha<1$, we obtain
 $$\left\|e^{it\sqrt{G}}(\sqrt{G})^{-\alpha(n-1)}
\chi_a(\sqrt{G})-e^{it\sqrt{G_0}}(\sqrt{G_0})^{-\alpha(n-1)}
\chi_a(\sqrt{G_0})\right\|_{L^{p'}\to L^p}$$ $$\le \int_0^1\left\|
\Phi(t;\theta)\right\|_{L^{p'}\to L^p}\theta^{\alpha(n-1)-1}d\theta 
\le C|t|^{-\alpha(n-1)/2}\int_0^1\theta^{-\alpha}d\theta
\le C|t|^{-\alpha(n-1)/2},\eqno{(4.5)}$$
which implies (1.3).\\

{\it Proof of Theorem 4.1.} By standard interpolation arguments, 
it suffices to prove (4.1) for $p=+\infty$, $p'=1$ and $p=p'=2$, and
(4.2) for $p=2$ and $p=+\infty$. The $L^2\to L^2$
estimates follow from Theorem 3.1. The estimate (4.2) with $p=+\infty$, 
$\alpha=1$, follows from the following

\begin{prop} For every $t\neq 0$, $0<h\le 1$, 
$0<\epsilon\ll 1$, we have
$$\left\|\Phi(t;h)\langle x\rangle^{-(n-1+\epsilon)/2}
\right\|_{L^2\to L^\infty}\le Ch^{1-n/2}|t|^{-(n-1)/2},\eqno{(4.6)}$$
with a constant $C>0$ independent of $t$ and $h$.
\end{prop}

{\it Proof.} Using (3.4) together with (2.3), (2.28) and (2.31), we obtain
$\forall f\in L^2$, $g\in L^1$,
$$\left|\left\langle \Phi(t;h)\langle x\rangle^{-(n-1+\epsilon)/2}f,g
\right\rangle\right|$$ $$\le Ch^2\left\|\Phi(t;h)\langle x
\rangle^{-(n-1+\epsilon)/2}f\right\|_{L^\infty}\|g\|_{L^1}
 +  Ch^{1-n/2}|t|^{-(n-1)/2}\|f\|_{L^2}\|g\|_{L^1}$$ 
  $$+Ch\int_0^t\left|\left\langle
\widetilde\varphi_1(h\sqrt{G_0})\sin\left((t-\tau)
\sqrt{G_0}\right)Ve^{i\tau\sqrt{G}}\varphi(h\sqrt{G}) 
\langle x\rangle^{-(n-1+\epsilon)/2}f,g\right\rangle\right|d\tau.
\eqno{(4.7)}$$
Denote by $A(t)$ the integral in the RHS of (4.7). In view of 
(2.4), (3.2) and (3.18), we have
$$|t|^{(n-1)/2}A(t)\le C\int_0^{t/2}\left\|\langle x
\rangle^{-(1+\epsilon)/2}e^{i\tau\sqrt{G}}\varphi(h\sqrt{G}) 
\langle x\rangle^{-(n-1+\epsilon)/2}f\right\|_{L^2}$$
 $$\times |t-\tau|^{(n-1)/2}\left\|\langle x
\rangle^{-(n+\epsilon)/2}\sin\left((t-\tau)\sqrt{G_0}\right)
\widetilde\varphi_1(h\sqrt{G_0})g\right\|_{L^2}d\tau$$
  $$+C\int_0^{t/2}|t-\tau|^{(n-1)/2}\left\|\langle x
\rangle^{-(n-1+\epsilon)/2}e^{i(t-\tau)\sqrt{G}}\varphi(h\sqrt{G}) 
\langle x\rangle^{-(n-1+\epsilon)/2}f\right\|_{L^2}$$
 $$\times \left\|\langle x
\rangle^{-1-\epsilon}\sin\left(\tau\sqrt{G_0}\right)
\widetilde\varphi_1(h\sqrt{G_0})g\right\|_{L^2}d\tau$$
 $$\le C\gamma\int_{-\infty}^{\infty}\left\|\langle x
\rangle^{-(1+\epsilon)/2}e^{i\tau\sqrt{G}}\varphi(h\sqrt{G}) 
\langle x\rangle^{-(n-1+\epsilon)/2}f\right\|_{L^2}^2d\tau$$
 $$+C\gamma^{-1}\int_{-\infty}^{\infty}|\tau|^{n-1}\left\|\langle x
\rangle^{-(n+\epsilon)/2}\sin\left(\tau\sqrt{G_0}\right)
\widetilde\varphi_1(h\sqrt{G_0})g\right\|_{L^2}^2d\tau$$
 $$+ C\gamma\int_{-\infty}^{\infty}
\langle\tau\rangle^{-1-\epsilon}|t-\tau|^{n-1}\left\|\langle x
\rangle^{-(n-1+\epsilon)/2}e^{i(t-\tau)\sqrt{G}}\varphi(h\sqrt{G}) 
\langle x\rangle^{-(n-1+\epsilon)/2}f\right\|_{L^2}^2d\tau$$
 $$+C\gamma^{-1}\int_{-\infty}^{\infty}\langle\tau\rangle^{1+\epsilon}
\left\|\langle x
\rangle^{-1-\epsilon}\sin\left(\tau\sqrt{G_0}\right)
\widetilde\varphi_1(h\sqrt{G_0})g\right\|_{L^2}^2d\tau$$
 $$\le C\gamma\|f\|_{L^2}^2+C\gamma^{-1}h^{-n}
\|g\|_{L^1}^2
=2Ch^{-n/2}\|f\|_{L^2}\|g\|_{L^1},\eqno{(4.8)}$$
if we take $\gamma=h^{-n/2}\|g\|_{L^1}/\|f\|_{L^2}$. 
By (4.7) and (4.8), we conclude
$$\left\|\Phi(t;h)\langle x
\rangle^{-(n-1+\epsilon)/2}f\right\|_{L^\infty}\le 
Ch^2\left\|\Phi(t;h)\langle x
\rangle^{-(n-1+\epsilon)/2}f\right\|_{L^\infty}$$ $$+
 Ch^{1-n/2}|t|^{-(n-1)/2}\|f\|_{L^2},\quad
\forall f\in L^2,\,0<h\le 1.\eqno{(4.9)}$$
Hence, there exists a constant $0<h_0<1$ such that for $0<h\le h_0$
we can absorbe the first term in the RHS of (4.9), thus obtaining (4.6)
in this case. Let now $h_0\le h\le 1$. Without loss of generality we
may suppose $h=1$. Then, the only term we need to estimate is
$$\left\|\left(\varphi_1(\sqrt{G})-\varphi_1(\sqrt{G_0})\right)
e^{it\sqrt{G}}\varphi(\sqrt{G})\langle x
\rangle^{-(n-1+\epsilon)/2}f\right\|_{L^\infty}.$$
By (2.34) (with $p=+\infty$, $h=1$), the operator
$$\left(\varphi_1(\sqrt{G})-\varphi_1(\sqrt{G_0})\right)\langle x
\rangle^\delta:L^2\to L^\infty$$
is bounded. Thus, the problem is reduced to estimating the norm
$$\left\|\langle x\rangle^{-\delta}
e^{it\sqrt{G}}\varphi(\sqrt{G})\langle x
\rangle^{-(n-1+\epsilon)/2}f\right\|_{L^2},$$
which in view of Theorem 3.4 is upper bounded by $C
\langle t\rangle^{-(n-1)/2}\left\|f\right\|_{L^2}$. This completes
the proof of (4.6) for all $0<h\le 1$.
\eproof

We will now establish the estimate
$$\left\|\Phi(t;h)\right\|_{L^1\to L^\infty}
\le Ch^{1-n}|t|^{-(n-1)/2},\quad \forall t\neq 0,
\,0<h\le 1.\eqno{(4.10)}$$ 
By (2.2) and (2.31), we obtain 
$$\left\|\Phi_1(t;h)f\right\|_{L^\infty}
\le Ch^2\left\|\Phi(t;h)f\right\|_{L^\infty}
+Ch^{-(n-3)/2}|t|^{-(n-1)/2}\|f\|_{L^1},\quad
\forall f\in L^1.\eqno{(4.11)}$$ 
On the other hand, we have 
  $$|t|^{(n-1)/2}\left|\left\langle \Phi_2(t;h)f,g
\right\rangle\right|$$ $$\le |t|^{(n-1)/2}\int_0^t\left|\left\langle
\widetilde\varphi_1(h\sqrt{G_0})\sin\left((t-\tau)
\sqrt{G_0}\right)Ve^{i\tau\sqrt{G}}\varphi(h\sqrt{G}) 
f,g\right\rangle\right|d\tau$$
 $$\le C\int_0^{t/2}\left\|\langle x
\rangle^{-(1+\epsilon)/2}e^{i\tau\sqrt{G}}\varphi(h\sqrt{G}) 
f\right\|_{L^2}$$
 $$\times |t-\tau|^{(n-1)/2}\left\|\langle x
\rangle^{-(n+\epsilon)/2}\sin\left((t-\tau)\sqrt{G_0}\right)
\widetilde\varphi_1(h\sqrt{G_0})g\right\|_{L^2}d\tau$$
 $$+C\int_0^{t/2}|t-\tau|^{(n-1)/2}\left\|\langle x
\rangle^{-(n+\epsilon)/2}e^{i(t-\tau)\sqrt{G_0}}\varphi(h\sqrt{G_0}) 
f\right\|_{L^2}$$
 $$\times \left\|\langle x
\rangle^{-(1+\epsilon)/2}\sin\left(\tau\sqrt{G_0}\right)
\widetilde\varphi_1(h\sqrt{G_0})g\right\|_{L^2}d\tau$$
 $$+C\int_0^{t/2}|t-\tau|^{(n-1)/2}\left\|\langle x
\rangle^{-(n-1+\epsilon)/2}\Phi(t-\tau;h)f\right\|_{L^2}$$
 $$\times \left\|\langle x
\rangle^{-1-\epsilon}\sin\left(\tau\sqrt{G_0}\right)
\widetilde\varphi_1(h\sqrt{G_0})g\right\|_{L^2}d\tau$$
 $$\le C\gamma\int_{-\infty}^{\infty}\left\|\langle x
\rangle^{-(1+\epsilon)/2}e^{i\tau\sqrt{G}}\varphi(h\sqrt{G}) 
f\right\|_{L^2}^2d\tau$$
 $$+C\gamma^{-1}\int_{-\infty}^{\infty}|\tau|^{n-1}\left\|\langle x
\rangle^{-(n+\epsilon)/2}\sin\left(\tau\sqrt{G_0}\right)
\widetilde\varphi_1(h\sqrt{G_0})g\right\|_{L^2}^2d\tau$$
 $$+ C\gamma\int_{-\infty}^{\infty}|\tau|^{n-1}\left\|\langle x
\rangle^{-(n+\epsilon)/2}e^{i\tau\sqrt{G_0}}\varphi(h\sqrt{G_0}) 
f\right\|_{L^2}^2d\tau$$
 $$+C\gamma^{-1}\int_{-\infty}^{\infty}\left\|\langle x
\rangle^{-(1+\epsilon)/2}\sin\left(\tau\sqrt{G_0}\right)
\widetilde\varphi_1(h\sqrt{G_0})g\right\|_{L^2}^2d\tau$$
$$+ C\gamma\int_{-\infty}^{\infty}
\langle\tau\rangle^{-1-\epsilon}|t-\tau|^{n-1}\left\|
\langle x\rangle^{-(n-1+\epsilon)/2}\Phi(t-\tau;h)
f\right\|_{L^2}^2d\tau$$
 $$+C\gamma^{-1}\int_{-\infty}^{\infty}\langle\tau\rangle^{1+\epsilon}
\left\|\langle x
\rangle^{-1-\epsilon}\sin\left(\tau\sqrt{G_0}\right)
\widetilde\varphi_1(h\sqrt{G_0})g\right\|_{L^2}^2d\tau.
\eqno{(4.12)}$$
Observe that by (3.2) and (2.33) we have
$$\int_{-\infty}^{\infty}\left\|\langle x
\rangle^{-(1+\epsilon)/2}e^{i\tau\sqrt{G}}\varphi(h\sqrt{G}) 
f\right\|_{L^2}^2d\tau\le C h^{-n}\|f\|_{L^1}^2.
\eqno{(4.13)}$$
By (2.4), (4.6), (4.12) and (4.13), we obtain
$$|t|^{(n-1)/2}\left|\left\langle \Phi_2(t;h)f,g
\right\rangle\right|\le C\gamma h^{-n}
\|f\|_{L^1}^2+C\gamma^{-1}h^{-n}\|g\|_{L^1}^2
=2Ch^{-n}\|f\|_{L^1}\|g\|_{L^1},$$
if we take $\gamma=\|g\|_{L^1}/\|f\|_{L^1}$, and hence
 $$\left\|\Phi_2(t;h)f\right\|_{L^\infty}\le 
 Ch^{-n}|t|^{-(n-1)/2}\|f\|_{L^1},\quad
\forall f\in L^1,\,0<h\le 1.\eqno{(4.14)}$$
By (3.4), (4.11) and (4.14), we conclude 
$$\left\|\Phi(t;h)f\right\|_{L^\infty}\le 
Ch^2\left\|\Phi(t;h)f\right\|_{L^\infty}+
 Ch^{1-n}|t|^{-(n-1)/2}\|f\|_{L^1},\quad
\forall f\in L^1,\,0<h\le 1.\eqno{(4.15)}$$
Hence, there exists a constant $0<h_0<1$ such that for $0<h\le h_0$
we can absorbe the first term in the RHS of (4.15), thus obtaining (4.10)
in this case. Let now $h_0\le h\le 1$. Without loss of generality we
may suppose $h=1$. Then, the only term we need to estimate is
$$\left\|\left(\varphi_1(\sqrt{G})-\varphi_1(\sqrt{G_0})\right)
e^{it\sqrt{G}}\varphi(\sqrt{G})f\right\|_{L^\infty}.$$
As above, using (2.34) (with $p=+\infty$, $h=1$), 
we reduce the problem to estimating the norm
$$\left\|\langle x\rangle^{-\delta}
e^{it\sqrt{G}}\varphi(\sqrt{G})f\right\|_{L^2},$$
which in view of Proposition 4.2 and (2.3) is upper bounded by $C
|t|^{-(n-1)/2}\left\|f\right\|_{L^1}$. This completes
the proof of (4.10) for all $0<h\le 1$.
\eproof

G. Vodev, Universit\'e de Nantes,
 D\'epartement de Math\'ematiques, UMR 6629 du CNRS,
 2, rue de la Houssini\`ere, BP 92208, 44332 Nantes Cedex 03, France

e-mail: georgi.vodev@math.univ-nantes.fr

\end{document}